\documentclass[a4paper,11pt]{article}
\usepackage{amsmath,amsfonts,amsthm}
\usepackage[a4paper,centering,hscale=0.75,vscale=0.8]{geometry}

\newtheorem{prop}{Proposition}[section]
\newtheorem{thm}[prop]{Theorem}
\newtheorem{lemma}[prop]{Lemma}
\newtheorem{coroll}[prop]{Corollary}
\newtheorem{defi}[prop]{Definition}
\theoremstyle{remark}
\newtheorem{rmk}[prop]{Remark}
\theoremstyle{definition}

\numberwithin{equation}{section}

\renewcommand{\P}{\mathbb{P}}
\newcommand{\E}{\mathbb{E}}
\newcommand{\F}{\mathfrak{F}}
\renewcommand{\H}{\mathcal{H}}

\renewcommand{\L}{\mathcal{L}}
\newcommand{\erre}{\mathbb{R}}

\renewcommand{\epsilon}{\varepsilon}
\newcommand{\ip}[2]{\langle#1,#2\rangle}

\newcommand{\ds}{\displaystyle}
\newcommand{\tr}{\mathop{\mathrm{Tr}}\nolimits}

\newsymbol\lesssim 132E

\title{Regular dependence on initial data for stochastic evolution
  equations with multiplicative Poisson noise}

\author{Carlo Marinelli \and Claudia Pr\'ev\^ot \and Michael R\"ockner}

\date{\normalsize 11 August 2008}

\begin{document}
\maketitle

\begin{abstract}
  We prove existence, uniqueness and Lipschitz dependence on the
  initial datum for mild solutions of stochastic partial differential
  equations with Lipschitz coefficients driven by Wiener and Poisson
  noise. Under additional assumptions, we prove G\^ateaux and
  Fr\'echet differentiability of solutions with respect to the
  initial datum. As an application, we obtain gradient estimates
  for the resolvent associated to the mild solution. Finally, we prove
  the strong Feller property of the associated semigroup.
\end{abstract}

\section{Introduction}
We shall consider the mild formulation of a stochastic PDE of the form
\begin{equation}
  \label{eq:musso}
du(t) = [Au(t) + f(t,u(t))]\,dt + B(t,u(t))\,dW(t)
+ \int_Z G(t,u(t),z)\,\bar\mu(dt,dz),
\qquad u(0)=x,
\end{equation}
where $W$ and $\bar\mu$ are a Wiener process and a compensated Poisson
measure, respectively, on a Hilbert space, thus including a large
class of equations driven by Hilbert space-valued L\'evy noise, thanks
to the L\'evy-It\^o decomposition theorem. Precise assumptions on the
data of (\ref{eq:musso}) are given in the next section.

The main contribution of this paper is global well-posedness (i.e.
existence, uniqueness, and regular dependence on the initial datum of
a mild solution on any time interval $[0,T]$, $T<\infty$) in spaces of
c\`adl\`ag predictable processes whose supremum (in time) has finite
$p$-th moments. While the $L_2$ result was fully settled by Kotelenez
\cite{Kote-Doob} about twenty-five years ago, the lack of an $L_p$
theory has been pointed out more recently in \cite{AWZ}. The new tool
allowing the development of such a theory is an infinite dimensional
Bichteler-Jacod inequality, which also holds for stochastic
convolutions. The $L_p$ existence result allows us to prove first and
second order Fr\'echet differentiability of the solution with respect
to the initial datum (for first order G\^ateaux differentiability the
$L_2$ theory is enough, see also \cite{AMR2}). Moreover, these
differentibility results are a key tool to prove that, as long as the
noise has a Brownian component, the semigroup associated to the SPDE
is regularizing, in particular, that it has the strong Feller
property.  An essential ingredient to obtain this result is a formula
of Bismut-Elworthy type, which only holds under non-degeneracy
assumptions on $B$. Finally, we also obtain gradient estimates on the
resolvent associated to the SPDE.

The issues considered in this paper are by now classical for
stochastic PDE with Wiener noise (see
e.g. \cite{cerrai-libro,DP-K,DZ92,FriKno} and references therein), but
comparable results do not seem to be available in the more general
jump case considered here. In fact, it is fair to say that the theory
of stochastic PDEs driven by jump noise is not yet fully developed,
even though recent years have witnessed a growing interest in the
area: let us just mention, without any claim of completeness, the
recent monograph \cite{PZ-libro}, where the semigroup approach is
discussed, \cite{LR-heat} for an analytic approach based on
generalized Mehler semigroups, as well as the earlier important
contributions \cite{Gyo-semimg,MP} for the variational approach.

Let us also mention that differentiability properties of the solution
of stochastic PDE play an essential role in the study of the
associated Kolmogorov equations. This direction of research, while
thoroughly pursued in the case of Wiener noise (see
e.g. \cite{cerrai-libro,DP-K,ZLNM}), is still in its infancy for equations
with jumps, and we hope that our results will be useful in this
respect.

The paper is organized as follows: in Section \ref{sec:mr} precise
assumptions on the SPDE (\ref{eq:musso}) are given, and the main
results on well-posedness and regular dependence on the initial datum
are stated. In Section \ref{sec:aux} we prove a Bichteler-Jacod
inequality for infinite dimensional stochastic integrals with respect
to Poisson random measures, and we extend it to corresponding
stochastic convolutions. This result is essential in order to obtain
$L_p$ well-posedness, and it could be interesting in its own right. We
also recall some results on the differentiability of implicit
functions in Banach spaces, on which the proofs of regular dependence
heavily rely.  Section \ref{sec:proofs} contains the proofs of the
well-posedness and differentiability results. In Section \ref{sec:ex}
we obtain an analytic consequence of these results, that is gradient
estimates for the resolvent associated to the (solution of the)
SPDE. Finally, in Section \ref{sec:Feller} we show that the semigroup
associated to the SPDE is strong Feller, if $B$ is not degenerate.

Finally, we would like to mention that a part of the results of this
paper has been announced in \cite{Knoche-CRAS} (based on
\cite{Knoche-diss}). There was, however, an error both in the
formulation and proof in what was called there Burkholder-Davis-Gundy
inequality for Poisson integrals, on which all subsequent results
depended. One point of this paper is to correct this error. The
corresponding inequality is contained in Proposition
\ref{prop:BDGconv} below. Then, as described above, among other
things we prove that all results announced in \cite{Knoche-CRAS} hold.

Let us conclude this introduction with some words about notation.  By
$a \lesssim b$ we mean that there exists a constant $N$ such that $a
\leq Nb$. To emphasize that the constant $N$ depends on a parameter
$p$ we shall write $N(p)$ and $a \lesssim_p b$. Generic constants,
which may change from line to line, are denoted by $N$.  Given two
separable Banach spaces $E$, $F$ we shall denote the space of linear
bounded operators from $E$ to $F$ by $\mathcal{L}(E,F)$. Similarly, if
$H$ and $K$ are Hilbert spaces, we shall denote the space of
trace-class and Hilbert-Schmidt operators from $K$ to $H$ by
$\mathcal{L}_1(K,H)$ and $\mathcal{L}_2(K,H)$, respectively.
$\mathcal{L}_1^+$ stands for the subset of $\mathcal{L}_1$ consisting
of all positive operators. We shall write $\mathcal{L}_1(H)$ in place
of $\mathcal{L}_1(H,H)$, and similarly for the other spaces.  Given a
self-adjoint operator $Q\in\mathcal{L}_1^+(K)$, we denote by
$\mathcal{L}_2^Q(K,H)$ the set of all (possibly unbounded) operators
$B:Q^{1/2}K\to H$ such that $BQ^{1/2}\in\mathcal{L}_2(K,H)$. The norms
in $\L_2(K,H)$ and $\L_2^Q(K,H)$ will be denoted by $|\cdot|_2$ and
$|\cdot|_Q$, without explicitly indicating the dependence on the
spaces $K$ and $H$.
Lebesgue measure is denoted by $\mathrm{Leb}$, without mentioning the
underlying space if no misunderstanding can arise.
Given a function $\phi:E \to F$, we set
\[
[\phi]_1 := \sup_{\substack{x,y \in E\\x\neq y}}
                \frac{|\phi(x)-\phi(y)|}{|x-y|},
\]
and we denote by $\partial \phi:E \times E \to F$ the map
$(x,y) \mapsto \partial_y \phi(x)$, where the directional derivative
$\partial_y \phi(x)$ is defined by
\[
\partial_y \phi(x) := \lim_{h\to 0} \mathcal{Q}^h_y\phi(x)
:= \lim_{h \to 0} \frac{\phi(x+hy)-\phi(x)}{h}.
\]
We shall also use the symbol $\partial \phi(x)$ to denote the
G\^ateaux derivative, so $\partial \phi(x) \in \mathcal{L}(E,F)$,
defined by $y \mapsto \partial_y \phi(x)$.  Analogously, given a
function $\phi:E_1 \times E_2 \to F$, where $E_1$, $E_2$ are
further Banach spaces, we define the following directional derivatives
\begin{align*}
  \partial_{1,y} \phi(x_1,x_2) &= \lim_{h\to 0} \mathcal{Q}_{1,y}^h\phi(x_1,x_2)
      := \lim_{h\to 0} \frac{\phi(x_1+hy,x_2)-\phi(x_1,x_2)}{h}\\
  \partial_{2,z} \phi(x_1,x_2) &= \lim_{h\to 0} \mathcal{Q}_{2,z}^h\phi(x_1,x_2)
      := \lim_{h\to 0} \frac{\phi(x_1,x_2+hz)-\phi(x_1,x_2)}{h},
\end{align*}
and the corresponding maps
\begin{align*}
\partial_1 \phi:E_1\times E_2\times E_1 \ni 
(x_1,x_2,y) &\mapsto \partial_{1,y}\phi(x_1,x_2) \in F,\\
\partial_2 \phi:E_1\times E_2\times E_2 \ni
(x_1,x_2,z) &\mapsto \partial_{2,z}\phi(x_1,x_2) \in F.
\end{align*}
Partial G\^ateaux derivatives are denoted by the same symbols.
Fr\'echet differentials are denoted by $D$, with subscripts if
necessary.  Moreover, in view of the canonical isomorphism between
$\mathcal{L}(E,\mathcal{L}(E,F))$ and $\mathcal{L}^{\otimes 2}(E,F)$,
the space of bilinear maps from $E$ to $F$, we can and will consider
$D^2\phi$ as a map from $E$ to $\mathcal{L}^{\otimes 2}(E,F)$. The
space of $k$ times continuously differentiable maps from $E$ to $F$
will be denoted by $C^k(E,F)$, and simply by $C^k(E)$ if $F=\erre$.

We shall occasionally use the following standard notation for
stochastic integrals with respect to semimartingales and random
measures:
\[
\phi\cdot M(t) := \int_0^t \phi(s)\,dM(s),
\qquad
\phi \star \mu(t) := \int_0^t\!\int \phi(s,z)\,\mu(ds,dz).
\]

\section{Main results}\label{sec:mr}
Let us begin with stating our precise assumptions on equation (\ref{eq:musso}).
We are given two real separable Hilbert spaces $H$, $K$ and a filtered
probability space $(\Omega,\mathcal{F},\mathbb{F},\P)$,
$\mathbb{F}=(\mathcal{F}_t)_{t\in[0,T]}$, on which a Wiener process
with covariance operator $Q \in \L_1^+(K)$ is defined. Moreover, we
are given a measure space $(Z,\mathcal{Z},m)$ and a Poisson measure
$\mu$ on $[0,T]\times Z$, independent of $W$, defined on the same
stochastic basis. The compensator (dual predictable projection) of
$\mu$ is $\mathrm{Leb}\otimes m$, and the compensated measure
$\bar{\mu}$ is $\bar\mu:=\mu-\mathrm{Leb}\otimes m$.

Denoting the predictable $\sigma$-field by $\mathcal{P}$, we shall
assume throughout the paper that the following assumptions are
satisfied:
\begin{itemize}
\item[(i)] $A$ is the generator of a strongly continuous semigroup on $H$;
\item[(ii)] $f: \Omega \times [0,T] \times H \to H$ and $B: \Omega
  \times [0,T] \times H \to \mathcal{L}_2^Q(K,H)$ are $\mathcal{P}
  \times \mathcal{B}(H)$-measurable functions;
\item[(iii)] $G: \Omega \times [0,T] \times H \times Z \to Z$ is a 
$\mathcal{P} \times \mathcal{B}(H) \times \mathcal{Z}$-measurable function;
\item[(iv)] $x$ is an $H$-valued $\mathcal{F}_0$-measurable random variable.
\end{itemize}
Further assumptions on the data of the problem will be specified when
needed. For simplicity, we shall suppress explicit dependence on
$\omega \in \Omega$ of all random elements, if no confusion can
arise. Let us also recall that, by (i), there exist $M$, $\sigma \geq
0$, such that $|e^{tA}| \leq M e^{\sigma t}$. We set, for future
reference, $M_T:=M e^{\sigma T}$.

The concept of solution we shall work with and the spaces where
solutions are sought are defined next.
\begin{defi}
  A predictable process $u:[0,T] \to H$ is a mild solution
  of (\ref{eq:musso}) if it satisfies
  \begin{align*}
    u(t) = &e^{tA}x + \int_0^t e^{(t-s)A}f(s,u(s))\,ds
    + \int_0^t e^{(t-s)A} B(s,u(s))\,dW(s)\\
    &+ \int_{(0,t]}\int_Z e^{(t-s)A} G(s,u(s),z)\,\bar\mu(ds,dz)
  \end{align*}
  $\P$-a.s. for all $t\in[0,T]$, where at the same time we assume that
  all integrals on the right-hand side exist.
\end{defi}
We shall also write $u(x)$ to emphasize the dependence on the initial
datum, and $u(t,x)$ will stand for the value of $u(x)$ at time
$t\in[0,T]$.

In the following, for simplicity of notation, we shall often write
$\int_0^t$ in place of $\int_{(0,t]}$.

\begin{defi}
  Let $p \geq 2$. We shall denote by $\mathcal{H}_p(T)$ and
  $\mathbb{H}_p(T)$ the spaces of all predictable processes
  $u:[0,T]\to H$ such that
  \[
  |[u]|_p := \big( \sup_{t\in[0,T]} \E|u(t)|^p \big)^{1/p} <\infty,
  \]
  and
  \[
  \|u\|_p := \big( \E\sup_{t\in[0,T]} |u(t)|^p \big)^{1/p} < \infty,
  \]
  respectively.
\end{defi}
\noindent For reasons that will become apparent later, we shall also
need to consider the same spaces endowed with the equivalent norms
\[
|[u]|_{p,\lambda} := \Big(\sup_{t\in[0,T]} \E|e^{-\lambda t}u(t)|^p\Big)^{1/p},
\qquad
\|u\|_{p,\lambda} := \Big(\E\sup_{t\in[0,T]} |e^{-\lambda t}u(t)|^p\Big)^{1/p},
\]
with $\lambda>0$.
We shall also use the notation $\mathbb{L}_p$ to denote the set
$L_p(\Omega,\mathcal{F},\P;H)$.

The following well-posedness result in $\H_2(T)$ is quite simple to
prove and it essentially relies only on the isometric formula for
stochastic integrals with respect to compensated Poisson measures (see
also \cite{AMR2,PZ-libro} for similar results).
\begin{thm}     \label{thm:mild1} 
  Assume that $x \in \mathbb{L}_2$, $e^{sA}B(t,x) \in
  \mathcal{L}^Q_2(K,H)$, $e^{sA}G(t,x,\cdot) \in L_2(Z,m)$ for all
  $(s,t,x)\in[0,T]^2\times H$, and there exist $h\in L_1([0,T])$ and $a
  \in H$ such that
  \begin{align}
    |e^{sA}(f(t,x) - f(t,y))|^2 + |e^{sA}(B(t,x) - B(t,y))|^2_Q&\nonumber\\
    + \int_Z |e^{sA}(G(t,x,z) - G(t,y,z))|^2\,m(dz)
    &\leq h(s)|x-y|^2,\label{eq:ii}\\
    \label{eq:i}
    |e^{sA} f(t,a)|^2 + |e^{sA} B(t,a)|^2_Q
      + \int_Z |e^{sA} G(t,a,z)|^2\,m(dz) &\leq h(s).
  \end{align}
  $\P$-a.s. for all $x$, $y \in H$ and $s$, $t \in [0,T]$.
  Then equation (\ref{eq:musso}) admits a unique mild solution in
  $\mathcal{H}_2(T)$. Moreover, the solution map $x \mapsto u(x)$ is
  Lipschitz from $\mathbb{L}_2$ to $\H_2(T)$.
\end{thm}
As briefly mentioned above, this theorem is stated and proved for its
simplicity, even though a more refined result holds true. In fact, one
can look for mild solutions of (\ref{eq:musso}) in the smaller (and
more regular) spaces $\mathbb{H}_p(T)$, obtaining also that solutions
have c\`adl\`ag paths. The price to pay is that one has to find suitable
estimates to replace the isometry of the stochastic integral. In order
to obtain such estimates, we need to assume that $A$ is
$\eta$-$m$-dissipative, i.e. that $A-\eta I$ is $m$-dissipative
for some $\eta \geq 0$ (this is equivalent to assuming that
$|e^{tA}| \leq e^{\eta t}$ for all $t\geq 0$, i.e. that the semigroup
generated by $A$ is of quasi-contraction type). On the other hand,
Theorem \ref{thm:mild1} above holds without the quasi-dissipativity
condition on $A$.

Our first main result is the following theorem, where the solution map
is defined from $\mathbb{L}_p$ to $\mathbb{H}_p(T)$.
\begin{thm}     \label{thm:mild2}
  Let $p \geq 2$. Assume that $A$ is $\eta$-$m$-dissipative, $x \in
  \mathbb{L}_p$, and there exist $h\in L_p([0,T])$ and $a \in H$ such
  that
  \begin{multline} \label{eq:ii'}
    |e^{sA}(f(t,x) - f(t,y))| + |B(s,x) - B(s,y)|_Q\\
    + \max\big( |G(s,x,\cdot)-G(s,y,\cdot)|_{L_2(Z,m)},
    |G(s,x,\cdot)-G(s,y,\cdot)|_{L_p(Z,m)} \big) \leq h(s) |x-y|,
  \end{multline}
  \begin{equation}
    \label{eq:i'}
    |e^{sA} f(t,a)|^2 + |B(s,a)|^2_Q
      + \int_Z \big( |G(s,a,z)|^2+|G(s,a,z)|^p \big)\,m(dz) \leq h(s)
  \end{equation}
  $\P$-a.s. for all $x$, $y \in H$ and $s$, $t \in [0,T]$.
  Then equation (\ref{eq:musso}) admits a unique c\`adl\`ag mild solution
  in $\mathbb{H}_p(T)$. Moreover, the solution map $x \mapsto u(x)$ is
  Lipschitz from $\mathbb{L}_p$ to $\mathbb{H}_p(T)$.
\end{thm}
\begin{rmk}
  Before we proceed to state other results, we would like to make the
  following remarks:
\begin{itemize}
\item[(i)] Much more general existence and uniqueness results in
  $\mathbb{H}_2(T)$ were proved by Kotelenez \cite{Kote-Doob}, where
  noise terms driven by general locally square integrable martingales
  are allowed, as well as locally Lipschitz coefficients with linear
  growth.
\item[(ii)] One could also consider equations driven by martingales
  with independent increments, ``embedding'' equations driven by
  compensated Poisson random measures, using the equivalence result of
  Gy\"ongy and Krylov \cite{KryGyo1}. This approach, however, even
  though very powerful, would be less transparent, and for this reason
  we prefer to work directly with equations driven by a Wiener process
  and a compensated Poisson measure. Let us also recall that, if one
  only wants to obtain results in $\mathbb{H}_2(T)$, then general
  stochastic martingale measures are also allowed, appealing to the
  results in \cite{Kote-Doob} and to the above mentioned procedure
  developed in \cite{KryGyo1}.
\item[(iii)] If the coefficients of (\ref{eq:musso}) are independent
  of $\omega \in \Omega$, then one can obtain the Markov property of
  solutions in a standard way, e.g. following the method of
  \cite[Sect. 2.9]{krylov} -- see also \cite{Gyo-semimg},
  \cite{PZ-libro}.
\item[(iv)] It is not difficult to prove that mild solutions are weak
  solutions (in the sense of PDEs), as it follows, roughly speaking,
  by a suitable stochastic version of Fubini's theorem. More details
  can be found e.g. in \cite[Sect.~9.3]{PZ-libro}.
\end{itemize}
\end{rmk}

Under the additional assumption that the coefficients $f$, $B$ and $G$
are G\^ateaux differentiable, we obtain that the solution map enjoys
the same property. For this to hold, the simpler $\mathcal{H}_2(T)$
well-posedness suffices. In particular, no quasi-$m$-dissipativity
assumption on $A$ is needed.

From here until the end of this section we assume, for simplicity
only, that $f$, $B$ and $G$ are deterministic maps that do not depend
on time. Given a Banach space $E$, we shall denote the space of
functions $\phi:Z \to E$ such that $\int_Z |\phi|_E^p\,dm<\infty$
by $L_p(Z,m;E)$.
\begin{thm}    \label{thm:gattino}
  Under the hypotheses of Theorem \ref{thm:mild1}, assume that
  \begin{itemize}
  \item[(i)] $f$ is G\^ateaux differentiable with $\partial f \in
    C(H\times H,H)$;
  \item[(ii)] $B$ is G\^ateaux differentiable and $\partial B \in
    C(H\times H,\L_2^Q(K,H))$.
  \item[(iii)] the map $x \mapsto G(x,z)$ is G\^ateaux differentiable
    for all $s\in ]0,T]$ and $z\in Z$, and
    \[
    x \mapsto e^{sA}\partial_{1,y} G(x,z) \in C(H,H)
    \]
    for all $s\in ]0,T]$, $y\in H$, and $z \in Z$;
  \item[(iv)] one has
    \[
    x \mapsto e^{sA} \partial_{1,y} G(x,\cdot) \in C(H,L_2(Z,m;H))
    \]
    for all $s\in ]0,T]$ and $y \in H$.
  \end{itemize}
  Then the solution map $x \mapsto u(x):\mathbb{L}_2 \to
  \mathcal{H}_2(T)$ is G\^ateaux differentiable and $\partial u: (x,y)
  \mapsto
  \partial_y u(x) \in C(\mathbb{L}_2 \times
  \mathbb{L}_2,\mathcal{H}_2(T))$ is the mild solution of
\begin{multline}     \label{eq:diff1}
dv(t) = Av(t)\,dt
+ \partial f(u(t,x)) v(t)\,dt
+ \partial B(u(t,x)) v(t)\,dW(t)\\
+ \int_Z \partial_1 G(u(t,x),z)v(t)\,\bar{\mu}(dt,dz),
\qquad v(0)=y.
\end{multline}
Moreover, one has
\[
|[ \partial_y u(x) ]|_2 \leq N |y|_{\mathbb{L}_2}
\]
for all $x$, $y\in\mathbb{L}_2$, where the constant $N$, which does
not depend on $x$ and $y$, is the Lipschitz constant of the solution
map $\mathbb{L}_2 \ni x \mapsto u(x) \in \mathcal{H}_2(T)$.
\end{thm}

On the other hand, in order to obtain Fr\'echet differentiability of
the solution map, the full $\mathbb{H}_p(T)$ well-posedness result is
needed. At this point we would like to stress that the following two
theorems cannot be proved, to the best of our knowledge, on the basis
of the already known $\mathbb{H}_2(T)$ well-posedness, even if one is
interested only in the Fr\'echet differentiability of the solution map
from $H$ to $\mathcal{H}_2(T)$.

\begin{thm}     \label{thm:frocetto}
  Let $q > p \geq 2$, and assume that the hypotheses of Theorem
  \ref{thm:mild2} are satisfied with $p$ replaced by $q$. Moreover,
  assume that
  \begin{itemize}
  \item[(i)] $f \in C^1(H,H)$ and $B \in C^1(H,\L(K,H))$;
  \item[(ii)] $x \mapsto e^{tA}DB(x) \in
    C\big(H,\L(H,\L_2^Q(K,H))\big)$ for all $t \in ]0,T]$;
  \item[(iii)] $x \mapsto G(x,z) \in C^1(H,H)$ for all $z \in Z$.
  \end{itemize}
  Then the solution map $\mathbb{L}_p \ni x \mapsto u(x) \in
  \mathbb{H}_p(T)$ is G\^ateaux differentiable and $(x,y) \mapsto
  \partial_y u(x) \in C(\mathbb{L}_p \times
  \mathbb{L}_p,\mathbb{H}_p(T))$ is the mild solution of
  (\ref{eq:diff1}) in $\mathbb{H}_p(T)$. Moreover, one has
  \[
    \| \partial_y u(x) \|_p \leq N |y|_{\mathbb{L}_p}
  \]
  for all $x$, $y \in \mathbb{L}_p$, where $N$ denotes the Lipschitz
  constant of $\mathbb{L}_p \ni x \mapsto u(x) \in
  \mathbb{H}_p(T)$. Finally, if $x \in \mathbb{L}_q$, then $x \mapsto
  u(x) \in C^1(\mathbb{L}_q,\mathbb{H}_p(T))$. In particular, $x
  \mapsto u(x) \in C^1(H,\mathbb{H}_p(T))$.
\end{thm}

\begin{thm}      \label{thm:frocione}
  Let $q>2p\geq 4$. Under the hypotheses of Theorem \ref{thm:mild2}
  with $p$ replaced by $q$, assume that
  \begin{itemize}
  \item[(i)] $f \in C^2(H,H)$, $B \in C^2(H,\L(K,H))$, and there exists
    $C_1>0$ such that
    \[
    |D^2f(x)|+|D^2B(x)| \leq C_1 \qquad \forall x \in H;
    \]
  \item[(ii)] the map $x \mapsto G(x,z): H \to H$ is twice Fr\'echet
    differentiable for all $z \in Z$ and
    \[
    x \mapsto D^2_1G(x,z) \in C\big(H,\L(H,\L(H))\big)
    \]
    for all $z \in Z$;
  \item[(iii)] there exists $h_1 \in L_p(Z,m) \cap L_2(Z,m)$ such that
    \[
    |D^2_1G(x,z)(y_1,y_2)| \leq h_1(z) |y_1| |y_2|
    \]
    for all $x$, $y_1$, $y_2 \in H$ and $z \in Z$;
  \item[(iv)] there exists $k \in L_2([0,T])$ such that
    \[
    \big|e^{tA}D^2B(x)(y,z)\big|_Q \leq k(t)|y|\,|z|.
    \]
  \end{itemize}
  Then the Fr\'echet derivative $Du: \mathbb{L}_q \to
  \L(\mathbb{L}_q,\mathbb{H}_p(T))$ is G\^ateaux differentiable.  Let
  $x$, $y_1$, $y_2 \in \mathbb{L}_q$, and $w:=[\partial
  Du(x)](y_1,y_2) \equiv [\partial^2 u(x)](y_1,y_2)$, $v_1=Du(x)y_1$,
  $v_2=Du(x)y_2$. Then $w$ is the mild solution of
  \begin{multline}     \label{eq:frocione}
  dw(t) = \big[ Aw(t) + Df(u(t))w(t) + D^2f(u(t))(v_1(t),v_2(t)) \big]\,dt\\
  + \int_Z \big[ D_1G(u(t),z)w(t) + D_1^2G(u(t),z)(v_1(t),v_2(t))
           \big]\,\bar{\mu}(dt,dz),
  \qquad w(0)=0.
  \end{multline}
  Moreover, there exists a constant $N=N(T,p,q)>0$ such that
  \[
  \| \partial Du(x)(y_1,y_2) \|_p \leq N |y_1|_{\mathbb{L}_q}
  |y_2|_{\mathbb{L}_q}
  \]
  for all $y_1$, $y_2 \in \mathbb{L}_q$.  Finally, if $q>4p\geq 8$,
  then
  \[
  x \mapsto u(x) \in C^2(\mathbb{L}_q,\mathbb{H}_p(T)).
  \]
  In particular, the solution map belongs to $C^2(H,\mathbb{H}_p(T))$.
\end{thm}

\section{Auxiliary results}\label{sec:aux}
\subsection{$L_p$ estimates for stochastic convolutions}
In order to prove Theorem \ref{thm:mild2} we need to establish a
maximal inequality for stochastic convolutions with respect to a
compensated Poisson measure, which may be of independent interest. For
related estimates (which hold only for stochastic integrals) in the
finite dimensional case see \cite{BGJ} and references therein, and for
the special case of stochastic integrals with respect to L\'evy
processes \cite{JKMP,ProTal-Euler}. Maximal inequalities for
stochastic convolutions can be found e.g. in
\cite{Ichi,Kote-sub,Kote-Doob}. None of the latter results, however,
seems to be useful to obtain the estimates we need to establish
well-posedness in $\mathbb{H}_p(T)$.

Let us begin with a Bichteler-Jacod inequality for Poisson integrals.
\begin{lemma}     \label{lm:31}
  Let $p \geq 2$. Assume that $g:[0,T]\times Z \to H$ is
  a predictable process such that the expectation on the right-hand side of
  (\ref{eq:BJ}) below is finite. Then one has
  \begin{equation} \label{eq:BJ} \E \Big| \sup_{t\leq T}
    \int_{(0,t]}\!\int_Z g(s,z)\,\bar\mu(ds,dz)\Big|^p \leq N \E \int_0^T
    \Big[ \int_Z |g(s,z)|^p\,m(dz) + \Big(\int_Z
    |g(s,z)|^2\,m(dz)\Big)^{p/2} \Big]\,ds,
  \end{equation}
  where $N=N(p,T)$, and $(p,T)\mapsto N$ is continuous.
\end{lemma}
\begin{proof}
Setting $\phi:H \to \erre$, $\phi(x)=|x|^p$, we have that $\phi$ is
twice Fr\'echet differentiable with derivatives
\[
\phi'(x): \eta \to p|x|^{p-2}\ip{x}{\eta}
\]
and
\[
\phi''(x): (\eta,\zeta) \mapsto
p(p-2)|x|^{p-4}\ip{x}{\eta}\ip{x}{\zeta} + p|x|^{p-2}\ip{\eta}{\zeta},
\qquad x \neq 0,
\]
$\phi''(0)=0$.
Let us set $X= g \star \bar{\mu}$. Then It\^o's formula (see
e.g. \cite{Met}) yields
\begin{align}
|X(t)|^p =&\; p\int_0^t \ip{|X(s-)|^{p-2}X(s-)}{dX(s)}\nonumber\\
&+ \sum_{s\leq t} \big( |X(s)|^p - |X(s-)|^p 
- p |X(s-)|^{p-2} \ip{X(s-)}{\Delta X(s)} \big) \label{eq:litorina}
\end{align}
$\P$-a.s. for all $t \leq T$, where, as usual, $\Delta X(s):=X(s)-X(s-)$.
Applying Taylor's formula to the function $\phi$ we obtain
\begin{align*}
&|X(s)|^p - |X(s-)|^p -p |X(s-)|^{p-2} \ip{X(s-)}{\Delta X(s)}\\
&\qquad = \frac12 p(p-2) \big| X(s-)+\xi\Delta X(s) \big|^{p-4}
  \big\langle X(s-)+\xi\Delta X(s),\Delta X(s) \big\rangle^2\\
&\qquad\quad + \frac12 p \big| X(s-)+\xi\Delta X(s) \big|^{p-2}
|\Delta X(s)|^2\\
&\qquad \leq \frac12 p(p-1) \big| X(s-)+\xi\Delta X(s) \big|^{p-2}
|\Delta X(s)|^2,
\end{align*}
where $\xi \equiv \xi(s) \in ]0,1[$ (see e.g. \cite[Thm. 4.18.1]{DeMarco2I}).
Since $|X(s-)+\xi\Delta X(s)| \leq |X(s-)| + |\Delta
X(s)|$, we also have
\[
\big| X(s-)+\xi\Delta X(s) \big|^{p-2} \lesssim_p |X(s-)|^{p-2}
+ |\Delta X(s)|^{p-2} \leq 
X^*(s-)^{p-2}
+ |\Delta X(s)|^{p-2},
\]
where $X^*(s):=\sup_{r\leq s} |X(r)|$.

Let us now assume, for the time being, that $X$ is bounded
$\P$-a.s.. Then the first term on the right-hand side of
(\ref{eq:litorina}) is a martingale with expectation zero, and we obtain
\[
\E |X(t)|^p \leq 
N(p) \E\sum_{s\leq t} \big(X^*(s-)^{p-2}|\Delta X(s)|^2 + |\Delta X(s)|^p
\big).
\]
Therefore, recalling that the compensator of $\mu$ is $m \otimes
\mathrm{Leb}$ and using Young's inequality
\[
ab \leq \frac{a^{\frac{p}{p-2}}}{\frac{p}{p-2}}
+ \frac{b^{p/2}}{p/2},
\]
we get
\begin{align*}
\E |X(t)|^p \lesssim_p\;& \E\int_0^t \big[X^*(s-)^{p-2}
|g(s,\cdot)|^2_{L_2(Z,m)}
+ |g(s,\cdot)|^p_{L_p(Z,m)}\big]\,ds\\
\lesssim_p\;& \E \int_0^t \big[X^*(s)^p + |g(s,\cdot)|^p_{L_2(Z,m)}
+ |g(s,\cdot)|^p_{L_p(Z,m)}\big]\,ds.
\end{align*}
Doob's inequality then yields
\[
\E X^*(t)^p \lesssim_p \E\int_0^t \big[X^*(s)^p + |g(s,\cdot)|^p_{L_2(Z,m)}
+ |g(s,\cdot)|^p_{L_p(Z,m)}\big]\,ds,
\]
hence, thanks to Gronwall's inequality, we obtain (\ref{eq:BJ}).

In order to remove the assumption that $X$ is bounded almost surely,
we shall proceed in two steps. Assume first that $|g(s,z)| \leq N$
a.s. for all $(s,z) \in [0,T] \times Z$, and define the stopping times
\[
\tau_n = \inf\{t \geq 0: \; |X(t)|>n\} \wedge T.
\]
We clearly have $\tau_n \to T$ a.s. as $n\to\infty$ because $X$ is
bounded on compact intervals. We also have
\[
|\Delta X(t)| \leq \sup_{(s,z)\in[0,T]\times Z} |g(s,z)| \leq N
\qquad\text{a.s.}
\]
hence, setting $Y_n = (1_{]0,\tau_n]} g) \star \bar{\mu} \equiv
X(t\wedge\tau_n)$, one easily sees that
\[
|Y_n(t)| \leq |Y_n(t-)| + \sup_{t\leq T} |\Delta X(t)| \leq n+N,
\]
and, by Fatou's lemma and passing to the limit,
\begin{align*}
  \E X^*(T)^p &\leq \liminf_{n\to\infty} \E Y_n^*(T)^p\\
  &\lesssim_p \liminf_{n\to\infty} \E
    \int_0^T \Big[\int_Z|1_{]0,\tau_n]} g(s,z)|^p\,m(dz)
       + \Big(\int_Z|1_{]0,\tau_n]} g(s,z)|^2\,m(dz)\Big)^{p/2}\Big]\,ds\\
  &\leq \E \int_0^T \Big[\int_Z|g(s,z)|^p\,m(dz)
       + \Big(\int_Z|g(s,z)|^2\,m(dz)\Big)^{p/2}\Big]\,ds,
\end{align*}
which proves the claim if $g$ is a.s. bounded. The general case can be
proved by setting
\[
g_n(s,z) :=
\begin{cases}
\ds g(s,z), & \text{if } |g(s,z)| \leq n\\[4pt]
\ds n\frac{g(s,z)}{|g(s,z)|}, & \text{if } |g(s,z)| > n,
\end{cases}
\]
and $X_n:=g_n \star \bar{\mu}$, from which it is easy to prove that,
by (\ref{eq:BJ}), $\{X_n\}_{n\in\mathbb{N}}$ is a Cauchy sequence in
$\mathbb{H}_p(T)$ and $X_n \to X$ in $\mathbb{H}_p(T)$, with $X=g
\star \bar{\mu}$.  Using again Fatou's lemma and recalling that
(\ref{eq:BJ}) holds for bounded integrands, we have
\begin{align*}
  \E X^*(T)^p &\leq \liminf_{n\to\infty} \E X_n^*(T)^p\\
  &\lesssim_p \liminf_{n\to\infty} \E
    \int_0^T \Big[\int_Z|g_n(s,z)|^p\,m(dz)
       + \Big(\int_Z|g_n(s,z)|^2\,m(dz)\Big)^{p/2}\Big]\,ds\\
  &\leq \E \int_0^T \Big[\int_Z|g(s,z)|^p\,m(dz)
       + \Big(\int_Z|g(s,z)|^2\,m(dz)\Big)^{p/2}\Big]\,ds,
\end{align*}
which concludes the proof.
\end{proof}
\begin{rmk}
  The ``usual'' Bichteler-Jacod inequality for L\'evy integrals (see
  e.g. \cite{ProTal-Euler}) follows immediately by (\ref{eq:BJ}) and
  the L\'evy-It\^o decomposition.  Moreover, the proof we gave is
  different from the ones in the literature and, apart from holding also
  in infinite dimensions, has the peculiarity of avoiding completely
  the use of the Burkholder-Davis-Gundy's inequality.
\end{rmk}

Inequality (\ref{eq:BJ}) can be extended also to stochastic
convolutions, even though in general they are not martingales.
\begin{prop}     \label{prop:BDGconv}
  Let $A$ be $m$-dissipative on $H$ and $g$ satisfies the hypotheses
  of Lemma \ref{lm:31}.  Then for all $p\in [2,\infty)$ there exists a
  constant $N$ such that
  \begin{multline}
    \label{eq:BDGconv}
    \E \sup_{t\leq T} \Big| 
    \int_0^t\int_Z e^{(t-s)A} g(s,z)\,\bar\mu(ds,dz)\Big|^p\\
    \leq
    N \E \int_0^T \Big[ \int_Z |g(s,z)|^p\,m(dz)
    + \Big(\int_Z |g(s,z)|^2\,m(dz)\Big)^{p/2} \Big]\,ds,
  \end{multline}
  where $N$ depends continuously on $p$ and $T$ only.
\end{prop}
\begin{proof}
  We shall follow the approach of \cite{HauSei}. In particular, by
  Sz.-Nagy's theorem on unitary dilations, there exists a Hilbert
  space $\bar{H}$, with $H$ isometrically embedded into $\bar{H}$, and a
  unitary strongly continuous group $T(t)$ on $\bar{H}$ such that $\pi
  T(t)x=e^{tA}x$ for all $x\in H$, $t\in \erre$, where $\pi$ denotes
  the orthogonal projection from $\bar{H}$ to $H$. Then we have,
  recalling that the operator norms of $\pi$ and $T(t)$ are less than or
  equal to one,
  \begin{eqnarray*}
  \E \sup_{t\leq T} 
  \Big| \int_0^t\int_Z e^{(t-s)A} g(s,z)\,\bar\mu(ds,dz)\Big|^p_H
  &=& 
  \E \sup_{t\leq T} \Big|
  \pi T(t) \int_0^t\int_Z T(-s) g(s,z)\,\bar\mu(ds,dz)\Big|^p_{\bar H}\\
  &\leq& |\pi|^p \; \sup_{t\leq T} |T(t)|^p \;
  \E \sup_{t\leq T} \Big|
  \int_0^t\int_Z T(-s) g(s,z)\,\bar\mu(ds,dz)\Big|^p_{\bar H}\\
  &\leq& \E \sup_{t\leq T} \Big|
  \int_0^t\int_Z T(-s) g(s,z)\,\bar\mu(ds,dz)\Big|^p_{\bar H}
  \end{eqnarray*}
  Since the integral in the last expression is a martingale,
  inequality (\ref{eq:BJ})
  implies that there exists a constant $N=N(p,T)$ such that
  \begin{eqnarray*}
  \lefteqn{\E \sup_{t\leq T} 
  \Big| \int_0^t\int_Z e^{(t-s)A} g(s,z)\,\bar{\mu}(ds,dz)\Big|^p}\\
  &\leq& N \E \int_0^T \Big[ \int_Z |T(-s)g(s,z)|^p \,m(dz)
    + \Big(\int_Z |T(-s)g(s,z)|^2 \,m(dz)\Big)^{p/2} \Big]\,ds\\
  &\leq& N \E \int_0^T \Big[ \int_Z |g(s,z)|^p\,m(dz)
    + \Big(\int_Z |g(s,z)|^2\,m(dz)\Big)^{p/2} \Big]\,ds
  \end{eqnarray*}
  where we have used again that $T(t)$ is a unitary group and that the
  norms of in $\bar{H}$ and $H$ are
  equal.
\end{proof}
\begin{coroll}
  Let $A$ be $\eta$-$m$-dissipative. Then inequality
  (\ref{eq:BDGconv}) holds, with $N$ a continuous function of $p$,
  $T$, and $\eta$.
\end{coroll}
\begin{proof}
  Follows by exactly the same arguments used above applied to the
  $m$-dissipative operator $A-\eta I$.
\end{proof}

\subsection{Differentiability of implicit functions}
In order to prove regular dependence of solutions with respect to the
initial datum, we shall need the following versions of the implicit
function theorem. Similar results can be found in the literature (see
e.g. \cite{cerrai-libro, DZ96, DZ02}), but we have included the
complete statements here for the reader's convenience. A proof of
these specific versions can be found in \cite{FriKno}.

Let $E$, $\Lambda$ be two Banach spaces, and $\Phi:\Lambda \times E
\to E$ a function such that
\[
|\Phi(\lambda,x)-\Phi(\lambda,y)| \leq \alpha |x-y|
\]
for all $\lambda \in \Lambda$ and all $x$, $y \in E$, with
$\alpha \in [0,1[$. Banach's fixed point
theorem implies the existence and uniqueness of a function
$\phi:\Lambda \to E$ such that $\Phi(\lambda,\phi(\lambda))=\phi(\lambda)$
for all $\lambda \in \Lambda$.
\begin{thm}     \label{thm:impl}
  Assume that $\lambda \mapsto \Phi(\lambda,x)$ is continuous for all
  $x\in E$. Then $\phi \in C(\Lambda,E)$. Moreover, if $\Phi$ is
  Lipschitz with respect to $\lambda$ uniformly over $x\in E$, then
  $\phi$ is Lipschitz.
\end{thm}
\begin{thm}     \label{thm:gatto-impl}
  Assume that $\Phi(\cdot,x):\Lambda \to E$ is continuous for all $x\in
  E$, and that the maps $\partial_1 \Phi: \Lambda \times E \times \Lambda
  \to E$, $\partial_2 \Phi: \Lambda \times E \times E \to E$ are
  continuous. Then $\phi$ is G\^ateaux differentiable and
  $(\lambda,\mu)\mapsto \partial_\mu\phi(\lambda)$ is continuous from
  $\Lambda \times \Lambda$ to $E$. Moreover, one has
\[
\partial_\mu\phi(\lambda) =
\big( I - \partial_2 \Phi(\lambda,\phi(\lambda))\big)^{-1}
\partial_{1,\mu} \Phi(\lambda,\phi(\lambda)).
\]
\end{thm}
In the formulation of the following theorems we shall denote by
$\Lambda_0$ and $E_0$ two Banach spaces continuously embedded in
$\Lambda$ and $E$, respectively. Moreover, $\Lambda_1$ will denote a
further Banach space continuously embedded in $\Lambda_0$.
\begin{thm}     \label{thm:frechet-impl}
  Assume that $\Phi$ satisfies the hypotheses of Theorem
  \ref{thm:gatto-impl}, also with $\Lambda_0$ and $E_0$ replacing
  $\Lambda$ and $E$, respectively. Moreover, assume that $\partial_1
  \Phi \in C(\Lambda_0 \times E_0,\L(\Lambda_0,E))$ and $\partial_2
  \Phi \in C(\Lambda_0 \times E_0,\L(E_0,E))$. Then $\partial\phi \in
  C(\Lambda_0,\L(\Lambda_0,E))$, hence $\phi \in C^1(\Lambda_0,E)$.
\end{thm}
\begin{thm}     \label{thm:gattone-impl} 
  Assume that both $\Phi:\Lambda \times E \to E$ and $\Phi:\Lambda_0
  \times E_0 \to E_0$ satisfy the hypotheses of Theorem
  \ref{thm:gatto-impl}. If $\Phi:\Lambda_0 \times E_0 \to E$ admits
  second-order directional derivatives, then $\phi:\Lambda_0 \to E$ is
  twice G\^ateaux differentiable with $\partial^2\phi \in
  C(\Lambda_0^3,E)$ and
\begin{align*}
\partial^2\phi(\lambda_0): (\mu_0,\nu_0) \mapsto&
\big( I - \partial_2 \Phi(\lambda_0,\phi(\lambda_0)) \big)^{-1}
\big[\partial_1^2\Phi(\lambda_0,\phi(\lambda_0))(\mu_0,\nu_0)\\
&\quad + \partial_1\partial_2\Phi(\lambda_0,\phi(\lambda_0))
                             (\partial_{\mu_0}\phi(\lambda_0),\nu_0)\\
&\quad + \partial_2\partial_1\Phi(\lambda_0,\phi(\lambda_0))
                             (\mu_0,\partial_{\nu_0}\phi(\lambda_0))\\
&\quad + \partial_2^2\Phi(\lambda_0,\phi(\lambda_0))
             (\partial_{\mu_0}\phi(\lambda_0),\partial_{\nu_0}\phi(\lambda_0))
 \big].
\end{align*}
\end{thm}
\begin{thm}     \label{thm:frocione-impl} 
  Assume that both $\Phi:\Lambda \times E \to E$ and $\Phi:\Lambda_0
  \times E_0 \to E_0$ satisfy the hypotheses of Theorem
  \ref{thm:gatto-impl}. Moreover, assume that $\Phi \in C^2(\Lambda_0
  \times E_0,E)$ and that $\phi \in C^1(\Lambda_1,E_0)$. Then the
  Fr\'echet derivative $D\phi:\Lambda_1 \to \mathcal{L}(\Lambda_1,E)$
  is G\^ateaux differentiable. Furthermore, if $\partial D\phi$ can be
  realized as a map $\Lambda_1 \to
  \mathcal{L}\big(\Lambda_1,\mathcal{L}(\Lambda_1,E_0)\big)$, then
  $\phi \in C^2(\Lambda_1,E)$.
\end{thm}
\begin{coroll}     \label{cor:frocione-impl}

  Let $\Phi$ be as in the previous theorem and $\phi \in
  C^1(\Lambda_1,E_0)$. Moreover, assume that $D\phi$ and $D_iD_j\Phi$,
  $i,j\in\{1,2\}$, are bounded. Then $\partial D\phi:\Lambda_1 \to
  \L(\Lambda_1,\L(\Lambda_1,E))$.
\end{coroll}

\subsection{Some regularization results}
We record for future reference some simple regularization and
approximation results which are used in the proofs of the main
results.
\begin{prop}     \label{prop:yosi}
  Let $u$ be the mild solution of (\ref{eq:musso}) in
  $\mathcal{H}_2(T)$, and $u_\lambda$ the strong solution of the
  equation
\begin{equation}     \label{eq:yosi}
du(t) = (A_\lambda u(t) + f(u(t))\,dt + B(u(t))\,dW(t)
+ \int_Z G(u(t-),z)\,\bar{\mu}(dt,dz),
\qquad u(0)=x,
\end{equation}
where $A_\lambda$ stands for the Yosida approximation of $A$. Then
$u_\lambda \to u$ in $\mathcal{H}_2(T)$ as $\lambda \to 0$.
\end{prop}
\begin{proof}
  We sketch the proof only, as it resembles the corresponding proof
  for equations driven by Wiener noise only.  In fact, the strong
  solution $u_\lambda$ of (\ref{eq:yosi}) is an adapted c\`adl\`ag
  process, and the predictable process $t \mapsto u_\lambda(t-)$ is a
  mild solution of (\ref{eq:yosi}). Recalling that, for a fixed
  $t\in[0,T]$, one has $u_\lambda(t)-u_\lambda(t-)=0$ almost surely
  (no jumps at a fixed time can occurr), we can proceed along the
  lines of e.g. \cite[Thm.~3.5]{DP-K}.
\end{proof}
\noindent In the following proposition we take $f$, $f_\varepsilon$, $B$,
$B_\varepsilon$, $G$, $G_\varepsilon$ independent of $t\in[0,T]$
and $\omega\in\Omega$.
\begin{prop}     \label{prop:vava}
  Let $u$ and $u_\varepsilon$ be, respectively, the mild solutions in
  $\mathcal{H}_2(T)$ of (\ref{eq:musso}) and of the equation obtained
  replacing $f$, $B$, and $G$ with $f_\varepsilon$, $B_\varepsilon$,
  and $G_\varepsilon$ in (\ref{eq:musso}), where $f_\varepsilon(x) \to f(x)$
  in $H$, $e^{tA}B_\varepsilon(x) \to e^{tA}B(x)$ in
  $\mathcal{L}^Q_2(H)$, and
  \[
  \int_Z \big| e^{tA}(G_\varepsilon(x,z)-G(x,z))\big|^2\,m(dz) \to 0
  \]
  for all $x \in H$ as $\varepsilon \to 0$.
  Moreover, assume that there exists $K > 0$ such that
  \[ 
  [f_\varepsilon]_1^2 +
  \big|e^{tA}(B_\varepsilon(x)-B_\varepsilon(y))\big|_Q^2 + \int_Z
  \big|e^{tA}(G_\varepsilon(x,z)-G_\varepsilon(y,z))\big|^2\,m(dz) \leq
  K|x-y|^2
  \]
  for all $t\in[0,T]$.  Then $u_\varepsilon \to u$ in
  $\mathcal{H}_2(T)$.
\end{prop}
\begin{proof}
  We have
  \begin{align*}
    \E|u_\varepsilon(t)-u(t)|^2 \;&\lesssim\; \E\int_0^t
    \big|e^{(t-s)A} [f_\varepsilon(u_\varepsilon(s)) - f(u(s))]
    \big|^2\,ds\\
    &\qquad + \E\int_0^t \big| e^{(t-s)A}
    [B_\varepsilon(u_\varepsilon(s)) - B(u(s))]
    \big|_Q^2\,ds\\
    &\qquad + \E\int_0^t\!\int_Z \big| e^{(t-s)A}
    [G_\varepsilon(u_\varepsilon(s),z)
    - G(u(s),z)] \big|^2\,m(dz)\,ds\\
    &=: I_1(\varepsilon) + I_2(\varepsilon) + I_3(\varepsilon),
  \end{align*}
  and
  \begin{align*}
  I_2(\varepsilon) &\lesssim \E\int_0^t \big| e^{(t-s)A}
  [B_\varepsilon(u_\varepsilon(s)) - B_\varepsilon(u(s))] \big|_Q^2\,ds
  + \E\int_0^t \big| e^{(t-s)A} [B_\varepsilon(u(s))
                   - B(u(s))] \big|_Q^2\,ds\\
  &\leq K \E\int_0^t |u_\varepsilon(s)-u(s)|^2\,ds + \delta_2(\varepsilon),
  \end{align*}
  where $\delta_2(\varepsilon) \to 0$ as $\varepsilon \to 0$, in view
  of the assumptions on $B_\varepsilon$ and by dominated convergence.
  Completely similar estimates can be obtained for $I_1(\varepsilon)$
  and $I_3(\varepsilon)$. We thus get
  \[
  \E|u_\varepsilon(t)-u(t)|^2 \leq N \int_0^t \E|u_\varepsilon(t)-u(t)|^2
  + \delta(\varepsilon),
  \]
  with $\delta(\varepsilon) \to 0$ as $\varepsilon \to 0$, and the
  conclusion follows by Gronwall's lemma.
\end{proof}


\section{Proofs}      \label{sec:proofs}
\begin{proof}[Proof of Theorem \ref{thm:mild1}]
  We sketch the proof only, as we follow the well-known approach based
  on Banach's fixed point theorem. We have to prove that the mapping
  $\F: \H_2(T) \to \H_2(T)$ defined by
\begin{eqnarray}
\F u(t) &=&  e^{tA} x 
+ \int_0^t e^{(t-s)A} f(s,u(s))\,ds
+ \int_0^t e^{(t-s)A} B(s,u(s))\,dW(s)\nonumber\\
&& + \int_0^t \int_Z e^{(t-s)A} G(s,u(s),z)\,\bar\mu(ds,dz)
\label{eq:fi}
\end{eqnarray}
is well defined and is a contraction, after which the result follows
easily. Let us show that, for any $u\in \H_2(T)$, $\F u$ admits a
predictable modification such that $|[\F u]|_2<\infty$. Predictability
of $\F u$ follows by the mean-square continuity of the stochastic
convolution term with respect to $\bar\mu$ in (\ref{eq:fi}): in fact,
setting $M_A(t)=\int_0^t\int_Z e^{(t-s)A}
G(s,u(s),z)\,\bar\mu(ds,dz)$, a simple calculation shows that, for
$0\leq s\leq t \leq T$,
\begin{eqnarray*}
\E|M_A(t)-M_A(s)|^2 &\lesssim&
\E \int_0^s\int_Z |e^{(t-r)A}-e^{(s-r)A}|^2|G(r,u(r),z)|^2\,m(dz)\,dr \\
&& + \E \int_s^t\int_Z |e^{(t-r)A}|^2 |G(r,u(r),z)|^2\,m(dz)\,dr,
\end{eqnarray*}
which converges to zero as $s\to t$.

\noindent Moreover, we have
\begin{multline*}
|[\F u]|^2_2 \;\lesssim\; \sup_{t\leq T} \E|e^{tA} x|^2
+ \sup_{t\leq T} \E\Big| \int_0^t e^{(t-s)A} f(s,u(s))\,ds \Big|^2\\
+ \sup_{t\leq T} \E\Big| \int_0^t e^{(t-s)A} B(s,u(s))\,dW(s) \Big|^2
+ \sup_{t\leq T} \E\Big| \int_0^t\int_Z 
      e^{(t-s)A} G(s,u(s),z)\,\bar\mu(ds,dz) \Big|^2.
\end{multline*}
Using the isometry for stochastic integrals, and noting that 
hypotheses (\ref{eq:i}) and (\ref{eq:ii}) imply the estimate
\[
\int_Z |e^{sA}G(t,x,z)|^2\,m(dz) \leq Nh(s)(1+|x|)^2,
\]
we have
\begin{eqnarray*}
\lefteqn{\sup_{t\leq T}\E\Big| \int_0^t\int_Z
e^{(t-s)A} G(s,u(s),z)\,\bar\mu(ds,dz) \Big|^2}\\
&=& \sup_{t\leq T} \E\int_0^t\int_Z |e^{(t-s)A} G(s,u(s),z)|^2\,m(dz)\,ds
\lesssim \sup_{t\leq T} \E\int_0^t h(t-s) (1+|u(s)|)^2\,ds\\
&\leq& 2|h|_{L_1} (1+\sup_{t\leq T} \E|u(t)|^2) < \infty.
\end{eqnarray*}
Analogous estimates for the remaining terms in (\ref{eq:fi}) are
classical (see e.g. \cite{DZ96}), hence $|[\F u]|_2^2 < \infty$.

We shall now prove that there exists $\lambda$ such that $|[\F u-\F
v]|_{2,\lambda} \leq N |[u-v]|_{2,\lambda}$, with $N<1$. In fact we
have
\begin{align*}
|[\F u -\F v]|^2_{2,\lambda} \;\lesssim\;&
  \sup_{t\leq T} e^{-2\lambda t}\E\Big| 
      \int_0^t e^{(t-s)A} \big(f(s,u(s)) - f(s,v(s))\big)\,ds \Big|^2\\
&+ \sup_{t\leq T} e^{-2\lambda t}\E\Big| 
      \int_0^t e^{(t-s)A} \big(B(s,u(s))-B(s,v(s))\big)\,dW(s) \Big|^2\\
&+ \sup_{t\leq T} e^{-2\lambda t}\E\Big| \int_0^t\int_Z e^{(t-s)A} 
         \big(G(s,u(s),z)-G(s,v(s),z)\big)\,\bar\mu(ds,dz) \Big|^2,
\end{align*}
and
\begin{align*}
  &\E\Big| \int_0^t\int_Z e^{(t-s)A} 
         \big(G(s,u(s),z)-G(s,v(s),z)\big)\,\bar\mu(ds,dz) \Big|^2\\
  &\qquad = \E \int_0^t\int_Z \big|e^{(t-s)A}
         \big(G(s,u(s),z)-G(s,v(s),z)\big)\big|^2\,m(dz)\,ds\\
  &\qquad \leq \E \int_0^t e^{2\lambda s} h(t-s) e^{-2\lambda s}|u(s)-v(s)|^2\,ds
     \leq |u-v|^2_{2,\lambda} \int_0^t e^{2\lambda s} h(t-s)\,ds\\
  &\qquad \leq e^{2\lambda t} |u-v|^2_{2,\lambda} \int_0^t
               e^{-2\lambda s} h(s)\,ds,
\end{align*}
which implies that the third summand on the right-hand side of the
previous estimate of $|[\F u -\F v]|^2_{2,\lambda}$ is bounded by
$|u-v|^2_{2,\lambda} \int_0^T e^{-2\lambda s}h(s)\,ds$, which
converges to zero as $\lambda \to \infty$.  Completely analogous
calculations for the other summands show that there exists
$N=N(T,h,\lambda)$ such that $|[\F u-\F v]|^2_{2,\lambda} \leq
N|[u-v]|^2_{2,\lambda}$, and that one can find $\lambda_0>0$ so that
$N(T,h,\lambda_0)<1$, thus obtaining, by Banach's fixed point theorem,
existence and uniqueness of a mild solution to (\ref{eq:musso}).
Finally, Lipschitz continuity of the solution map follows
by Theorem \ref{thm:impl}.
\end{proof}
\begin{rmk}     \label{rmk:lipco}
  One can also prove by a direct calculation that $x \mapsto u(x)$ is
  Lipschitz. This method has the advantage of yielding explicit
  estimates on the Lipschitz constant, and will be useful to establish
  the strong Feller property. In fact, one has
  \begin{eqnarray*}
  u(t,x)-u(t,y) &=& e^{tA}(x-y)
  + \int_0^t e^{(t-s)A} [f(s,u(s,x))-f(s,u(s,y))]\,ds\\
  && + \int_0^t e^{(t-s)A} [B(s,u(s,x))-B(s,u(s,y))]\,dW(s)\\
  && + \int_0^t\int_Z e^{(t-s)A} [G(s,u(s,x),z)-G(s,u(s,y),z)]\,\bar\mu(ds,dz),
  \end{eqnarray*}
  hence, squaring both sides and taking expectations,
  \[
  \E|u(t,x)-u(t,y)|^2 \leq 2 M^2 e^{2\eta t} \E|x-y|^2 +
  (2t+1)\int_0^t h(t-s)\E|u(s,x)-u(s,y)|^2\,ds,
  \]
  which yields, via Gronwall's inequality,
  \[
  |[u(x)-u(y)]|_2 \leq \sqrt{2} M e^{(\eta + |h|_1)T + |h|_1/2}
  |x-y|_{\mathbb{L}_2}.
  \]
\end{rmk}


\begin{proof}[Proof of Theorem \ref{thm:mild2}]
  We shall use a fixed point argument in the space $\mathbb{H}_p(T)$.
  In particular, we want to prove that the mapping $\F$ defined as in
  (\ref{eq:fi}) is a well-defined contraction on $\mathbb{H}_p(T)$.
  Here we limit ourselves to prove that there exists $N<1$ such that
  $\| \F u - \F v \|_{p,\lambda} \leq N\|u-v\|_{p,\lambda}$ for all
  $u$, $v \in \mathbb{H}_p(T)$, with a suitably chosen $\lambda \geq
  0$. In fact, this implies
  \[
  \|\F u\|_p \lesssim \|u - a\|_p + \|\F a\|_p < \infty
  \]
  for all $u \in \mathbb{H}_p(T)$, thanks to (\ref{eq:i'}). Moreover,
  predictability of $\F u$, $u\in \mathbb{H}_p(T)$, follows as in the
  proof of Theorem \ref{thm:mild1}.

  We have
  \begin{eqnarray*}
    \|\F u - \F v\|^p_{p,\lambda} &\lesssim_p& \E\sup_{t\leq T}\Big|
          e^{-\lambda t} \int_0^t e^{(t-s)A} (f(s,u(s))-f(s,v(s)))\,ds \Big|^p\\
    && + \E\sup_{t\leq T}\Big| e^{-\lambda t} \int_0^t
          e^{(t-s)A} (B(s,u(s))-B(s,v(s)))\,dW(s) \Big|^p\\
    && + \E\sup_{t\leq T}\Big| e^{-\lambda t} \int_0^t\int_Z
         e^{(t-s)A} (G(s,u(s),z)-G(s,v(s),z))\,\bar\mu(ds,dz) \Big|^p\\
    && =: A_1 + A_2 + A_3.
  \end{eqnarray*}
  The term $A_1$ on the right-hand side is bounded from above, thanks
  to (\ref{eq:ii'}) and H\"older's inequality, by
  \begin{align*}
    &T^{p-1} \E\sup_{t\leq T} e^{-p\lambda t} \int_0^t h(t-s)^p e^{p\lambda s}
     \big( e^{-\lambda s} |u(s)-v(s)|\big)^p\,ds\\
    &\quad \leq T^{p-1} \|u-v\|^p_{p,\lambda} \sup_{t\leq T}
           \int_0^t h(t-s)^p e^{-p\lambda(t-s)}\,ds\\
    &\quad \leq T^{p-1} |h_\lambda|^p_{L_p([0,T])} \|u-v\|^p_{p,\lambda}, 
  \end{align*}
  where $h_\lambda(s):=e^{-\lambda s}h(s)$.

  Moreover, since
  \[
    A_3 = \E\sup_{t\leq T} \Big| \int_0^t\int_Z e^{(t-s)(A-\lambda I)}
              e^{-\lambda s}(G(s,u(s),z)-G(s,v(s),z))\,\bar\mu(ds,dz) \Big|^p
  \]
  and, by a slight modification of the proof of (\ref{eq:BDGconv}),
  \begin{align*}
  &\E\sup_{t\leq T} \Big| \int_0^t\int_Z e^{(t-s)(A-\lambda I)}
              \phi(s,z)\,\bar\mu(ds,dz) \Big|^p\\
  &\quad \lesssim_{p,T,\eta} e^{-\lambda p T} \E\int_0^T e^{\lambda p s}
         \big( |\phi(s,\cdot)|^p_{L_p(Z,m)} + |\phi(s,\cdot)|^p_{L_2(Z,m)}
         \big)\,ds,
  \end{align*}
  we obtain, thanks to (\ref{eq:ii'}),
  \begin{align*}
  A_3 &\lesssim_{p,T,\eta} e^{-\lambda p T} \E\Big[\int_0^T e^{\lambda p s}
          e^{-\lambda p s}|G(s,u(s),\cdot)-G(s,v(s),\cdot)|^p_{L_2(Z,m)}\,ds\\
      &\quad\qquad + \int_0^T e^{\lambda p s} e^{-\lambda p s}
             |G(s,u(s),\cdot)-G(s,v(s),\cdot)|^p_{L_p(Z,m)} \,ds \Big]\\
      &\leq \|u-v\|^p_{p,\lambda} e^{-\lambda p T}
            \int_0^T e^{\lambda p s} h(s)^p\,ds
       \leq |\tilde{h}_\lambda|^p_{L_p([0,T])} \|u-v\|^p_{p,\lambda},
  \end{align*}
  where $\tilde{h}_\lambda(s):=e^{-\lambda s}h(T-s)$. Classical
  maximal inequalities for stochastic convolutions with respect to
  Wiener processes yield a completely analogous estimate for $A_2$.
  Observing that the norms of $h_\lambda$ and $\tilde{h}_\lambda$
  appearing in the above estimates tend to zero as $\lambda \to
  \infty$, we conclude that there exists a constant
  $N=N(\lambda,T,p,\eta)$ such that $\|\F u - \F v\|_{p,\lambda}
  \leq N\|u-v\|_{p,\lambda}$, with $N<1$ for some $\lambda>0$
  sufficiently large. The existence and uniqueness of a solution, as
  well as its Lipschitz continuity with respect to the initial datum,
  follows then by Banach's fixed point theorem, as for Theorem
  \ref{thm:mild1}, and by the equivalence of the norms
  $\|\cdot\|_{p,\lambda}$ for $\lambda \geq 0$.
\end{proof}
\begin{rmk}
  The Lipschitz continuity of the solution map, in analogy to the
  previous remark, could also be obtained by a direct calculation.
  However, in this case the norm of $\mathbb{H}_p(T)$ is somewhat more
  difficult to work with. In Section \ref{sec:ex} below we shall
  obtain some estimates for the Lipschitz constant of the solution map
  under additional assumptions on the coefficient of the SPDE.
\end{rmk}


\begin{proof}[Proof of Theorem \ref{thm:gattino}]
  It is enough to prove the statements in the case $B \equiv 0$. We
  are going to apply Theorem \ref{thm:gatto-impl}, with
  $\Lambda=\mathbb{L}_2$, $E=\mathcal{H}_2(T)$. The latter space needs
  to be endowed with a norm $|[\,\cdot\,]|_{2,\lambda}$, where $\lambda>0$
  is chosen in such a way that $\F:\mathbb{L}_2 \times
  \mathcal{H}_2(T) \to \mathcal{H}_2(T)$ is a contraction in its
  second argument. However, in view of the equivalence of the norms
  $|[\,\cdot\,]|_{2,\lambda}$, we shall perform the calculations assuming
  $\lambda=0$, without loss of generality.

  It is immediate that the directional derivative $\partial_{1,y}
  \mathfrak{F}(x,u)$ coincides with the map $t \mapsto e^{tA}y$, which
  clearly belongs to $\mathcal{H}_2(T)$.
Moreover, we have
\begin{multline*}
\Big|\!\Big[ \mathcal{Q}^h_{2,v} \F(x,u)
- \int_0^\cdot e^{(\cdot-s)A} \partial_{v(s)} f(u(s))\,ds
- \int_0^\cdot\!\int_Z e^{(\cdot-s)A} \partial_{1,v(s)}
G(u(s),z)\,\bar{\mu}(ds,dz) \Big]\!\Big|_2\\
\leq \Big|\!\Big[ \int_0^\cdot e^{(\cdot-s)A} \big[ \mathcal{Q}^h_{v(s)} f(u(s))
- \partial_{v(s)} f(u(s))\big]\,ds \Big]\!\Big|_2\\
+ \Big|\!\Big[ \int_0^\cdot\!\int_Z e^{(\cdot-s)A}
  \big[ \mathcal{Q}^h_{1,v(s)} G(u(s),z)
    - \partial_{1,v(s)} G(u(s),z)\big]\,\bar{\mu}(ds,dz) \Big]\!\Big|_2
\end{multline*}
The first term on the right hand side of this inequality tends to zero
as $h\to 0$ by obvious estimates and the dominated convergence
theorem. Using the isometric property of the stochastic integral, the
square of the second term is equal to
\[
\E \int_0^T \int_Z
\big|e^{(T-s)A}\mathcal{Q}^h_{1,v(s)}G(u(s),z) - 
e^{(T-s)A}\partial_{1,v(s)} G(u(s),z)\big|^2
\,m(dz)\,ds.
\]
In view of assumptions (ii)--(iv), a simple computation shows that
\[
\frac{e^{tA}G(x+hy,\cdot) - e^{tA}G(x,\cdot)}{h}
\to e^{tA}\partial_{1,y}G(x,\cdot) \qquad \forall x \in H
\]
in $L_2(Z,m;H)$ as $h \to 0$, whence we obtain convergence to zero of
the second term in the above estimate again by dominated convergence.
Therefore we have
\[
[\partial_{2,v} \F(x,u)](t) =
\int_0^t e^{(t-s)A} \partial_{v(s)} f(u(s))\,ds
+ \int_0^t\int_Z e^{(t-s)A} \partial_{1,v(s)} G(u(s),z)\,\bar{\mu}(ds,dz).
\]
The continuity of $\partial_1\F$ and $\partial_2\F$, considered as
maps $\mathbb{L}_2\times \mathcal{H}_2(T) \times \mathbb{L}_2 \to
\mathcal{H}_2(T)$ and $\mathbb{L}_2\times \mathcal{H}_2(T) \times
\mathcal{H}_2(T) \to \mathcal{H}_2(T)$, respectively, can be proved in a
completely similar way, and we omit the details.

Let us now prove the second assertion of the theorem: let $x$,
$y$ in $\mathbb{L}_2$. Then Theorem \ref{thm:gatto-impl} yields
\[
\partial_y u(x) = \big(I-\partial_2\F(x,u(x))\big)^{-1}
\partial_{1,y} \F(x,u(x)),
\]
thus also
\[
\partial_y u(x) = \partial_{1,y} \F(x,u(x))
+ \partial_2\F(x,u(x))\partial_y u(x),
\]
and the result follows substituting in the previous formula the
expressions for the directional derivatives of $\F$ found above.

The last assertion of the theorem is a direct consequence of the
definition of directional derivative and the fact that the solution
map $x \mapsto u(x)$ is Lipschitz.
\end{proof}


\begin{proof}[Proof of Theorem \ref{thm:frocetto}]
  One can prove that the solution map is G\^ateaux differentiable as
  in the proof of Theorem \ref{thm:gattino}, except for the fact that
  one cannot use the isometric property of the stochastic integral,
  but has to rely on the estimate (\ref{eq:BDGconv}).
  In particular, one has
  \begin{align*}
  &\Big\| \int_0^\cdot\int_Z e^{(t-s)A} \big[ \mathcal{Q}^h_{1,v(s)} G(u(s),z)
  - \partial_{1,v(s)} G(u(s),z)\big]\,\bar{\mu}(ds,dz) \Big\|^p_p\\
  &\qquad \lesssim \E\int_0^T \big| \mathcal{Q}^h_{1,v(s)} G(u(s),\cdot)
  - \partial_{1,v(s)} G(u(s),\cdot)\big|^p_{L_p(Z,m)}\,ds\\
  &\qquad\quad + \E\int_0^T \big| \mathcal{Q}^h_{1,v(s)} G(u(s),\cdot)
  - \partial_{1,v(s)} G(u(s),\cdot)\big|^p_{L_2(Z,m)}\,ds,
  \end{align*}
  which converges to zero by dominated convergence, thanks to the
  assumptions on $G$.

  In order to prove that $x \mapsto u(x)$ is also Fr\'echet
  differentiable, let us set $\Lambda_0=\mathbb{L}_q$,
  $\Lambda=\mathbb{L}_p$, $E_0=\mathbb{H}_q(T)$, and
  $E=\mathbb{H}_p(T)$, and apply Theorem \ref{thm:frechet-impl}.
  We are going to prove that the partial G\^ateaux derivatives
  \begin{align*}
    \partial_1\F: \mathbb{L}_q \times \mathbb{H}_q(T) &\to
    \mathcal{L}\big( \mathbb{L}_q,\mathbb{H}_p(T) \big),\\
    \partial_2\F: \mathbb{L}_q \times \mathbb{H}_q(T) &\to
    \mathcal{L}\big( \mathbb{H}_q(T),\mathbb{H}_p(T) \big)
\end{align*}
are continuous, which implies that $\F$ is Frech\'et differentiable
(see e.g. \cite{AmbPro}).
Since $\partial_1\F:(x,u) \mapsto e^{tA}x$ is clearly continuous, it
suffices to show that $\partial_2\F$ is continuous.
To this purpose, let $\{x_n\} \subset \mathbb{L}_q$,
$x \in \mathbb{L}_q$, $\{u_n\} \subset \mathbb{H}_q(T)$,
$u$, $w \in \mathbb{H}_q(T)$ such that $x_n \to x$ in
$\mathbb{L}_q$, $u_n \to u$ in $\mathbb{H}_p(T)$, and
$\|w\|_q \leq 1$.
We shall prove that
\[
\big\|
\partial_2 \F(x_n,u_n)w - \partial_2 \F(x,u)w
\big\|_p \to 0
\]
as $n \to \infty$. In fact, the above expression is no greater than
\begin{multline*}
\Big\|
\int_0^\cdot e^{(\cdot-s)A} \big[Df(u_n(s))w(s) - Df(u(s))w(s)\big]\,ds
\Big\|_p\\
+ \Big\|
\int_0^\cdot\int_Z e^{(\cdot-s)A} \big[D_1G(u_n(s),z)w(s)
    - D_1G(u(s),z)w(s)\big]\,\bar{\mu}(ds,dz) \Big\|_p =: I_1(n) + I_2(n),
\end{multline*}
and, using H\"older's inequality with conjugate exponents $r=q/p$ and $r'$,
\begin{align*}
I_1(n)^p &\lesssim \E\int_0^T \big| [Df(u_n(s)) - Df(u(s))] w(s)
\big|^p\,ds\\
&\leq \left( \E\int_0^T \big| Df(u_n(s)) - Df(u(s)) \big|^{pr'}\,ds
\right)^{1/r'}
\left( \E\int_0^T |w(s)|^{pr}\,ds\right)^{1/r}\\
&\lesssim \|w\|_q^p \left( \E\int_0^T \big| Df(u_n(s)) - Df(u(s))
\big|^{\frac{pq}{q-p}}\,ds \right)^{\frac{q-p}{q}}
\end{align*}
Since $|Df(u_n(s))-Df(u(s))| \leq 2[f]_1$, the dominated convergence theorem
and the continuity of $Df$ imply that $I_1(n) \to 0$ as $n \to \infty$.

\noindent Similarly, using the maximal inequality (\ref{eq:BDGconv}), we obtain
\begin{align*}
I_2(n)^p &\leq
\E \bigg| \sup_{t \leq T} \int_0^t\int_Z
e^{(t-s)A} \big[D_1G(u_n(s),z)w(s) - D_1G(u(s),z)w(s)\big]\,\bar{\mu}(ds,dz)
   \bigg|^p\\
&\lesssim
\E \int_0^T \big(\big|D_1G(u_n(s),\cdot) - D_1G(u(s),\cdot)\big|^p_{L_p(Z,m)}\\
&\qquad\qquad + \big|D_1G(u_n(s),\cdot) - D_1G(u(s),\cdot)\big|^p_{L_2(Z,m)}\big)
          |w(s)|^p \, ds,
\end{align*}
which converges to zero as $n \to \infty$ by arguments completely
analogous to the above ones.
\end{proof}

\begin{proof}[Proof of Theorem \ref{thm:frocione}]
  We are going to apply Theorem \ref{thm:frocione-impl}, with
  $\Lambda=\mathbb{L}_p$, $\Lambda_0=\mathbb{L}_{q'}$,
  $\Lambda_1=\mathbb{L}_q$, and $E=\mathbb{H}_p(T)$,
  $E_0=\mathbb{H}_{q'}(T)$. Here $q' \in ]2p,q[$. In analogy to what
  we have done before, we shall endow $E$ and $E_0$ with norms
  $\|\cdot\|_{q',\lambda(q')}$ and $\|\cdot\|_{p,\lambda(p)}$,
  respectively, where the $\lambda$ are chosen in a such a way that
  $\F:\mathbb{L}_r \times \mathbb{H}_r(T) \to \mathbb{H}_r(T)$ are
  contractions in the second argument, but we shall perform the
  calculations assuming $\lambda \equiv 0$, without loss of generality.

  It is clear that, in view of Theorem \ref{thm:frocetto}, it is
  enough to prove that $\F \in C^2(\mathbb{L}_{q'}\times
  \mathbb{H}_{q'}(T),H_{p}(T))$.
  Since $\partial_1D_1\F \equiv \partial_1D_2\F \equiv \partial_2D_1\F
  \equiv 0$, we only have to consider $\partial_2 D_2 \F$. Let us
  first prove that
  \begin{align}
    [\partial_{2,v} D_2\F(x,u)w](t) =&
    \int_0^t e^{(t-s)A} D^2f(u(s))(v(s),w(s))\,ds\nonumber\\
    &+ \int_0^t\!\int_Z e^{(t-s)A} D_1^2 G(u(s),z)(v(s),w(s))
          \,\bar{\mu}(ds,dz), \qquad t \in [0,T],\label{eq:dido}
  \end{align}
  for all $v$, $w \in \mathbb{H}_{q'}(T)$. In fact, the $p$-th power
  of the $\mathbb{H}_p(T)$ norm of the difference between
  $\mathcal{Q}_{2,v}^h D_2 \F(x,u)w$ and the right-hand side of
  (\ref{eq:dido}) is not greater than a constant times
  \begin{align*}
    \E &\sup_{t\leq T} \Big| \int_0^t e^{(t-s)A} \big[
       \mathcal{Q}^h_{v(s)} Df(u(s))w(s) - D^2f(u(s))(v(s),w(s))\big]\,ds
                     \Big|^p\\
    &+ \E \sup_{t\leq T} \Big| \int_0^t\!\int_Z e^{(t-s)A} \big[
       \mathcal{Q}^h_{1,v(s)} D_1G(u(s),z)w(s) - D_1^2G(u(s),z)(v(s),w(s))\big]
              \,\bar{\mu}(ds,dz) \Big|^p\\
    &=: I_1(h) + I_2(h).
  \end{align*}
  We have, by H\"older's inequality,
  \begin{align*}
    I_1(h) &\lesssim \E \int_0^T \big|\big( \mathcal{Q}^h_{v(s)} Df(u(s))
    - D^2f(u(s))(v(s),\cdot)
    \big)w(s) \big|^p\,ds\\
    &\lesssim \|w\|_{q'}^{p} \Big( \E \int_0^T \big|
    \mathcal{Q}^h_{v(s)} Df(u(s)) - D^2f(u(s))(v(s),\cdot)
    \big|^{2p}_{\mathcal{L}(H)}\,ds \Big)^{1/2},
  \end{align*}
  which converges to zero as $h \to 0$ by dominated convergence,
  thanks to the boundedness of $D^2f$.
  Similarly, using inequality (\ref{eq:BDGconv}), we get
  \begin{align*}
    I_2(h) &\lesssim \E \int_0^T \big| \mathcal{Q}^h_{1,v(s)} D_1G(u(s),\cdot)w(s)
       - D_1^2G(u(s),\cdot)(v(s),w(s)) \big|^p_{L_p(Z,m)}\,ds\\
        &\qquad + \E \int_0^T \big| \mathcal{Q}^h_{1,v(s)} D_1G(u(s),\cdot)w(s)
       - D_1^2G(u(s),\cdot)(v(s),w(s)) \big|^p_{L_2(Z,m)}\,ds\\
        &= I_{21}(h) + I_{22}(h),
  \end{align*}
  and, by H\"older's inequality,
  \begin{align*}
    I_{21}(h) &\leq \E\int_0^T \big| \mathcal{Q}^h_{1,v(s)}
    D_1G(u(s),\cdot) - D_1^2G(u(s),\cdot)(v(s),\cdot)
    \big|^p_{L_p(Z,m;\mathcal{L}(H))} ds
    |w(s)|^p\,ds\\
    &\lesssim \|w\|_{q'}^p \, \E\int_0^T \big| \mathcal{Q}^h_{1,v(s)}
    D_1G(u(s),\cdot) - D_1^2G(u(s),\cdot)(v(s),\cdot)
    \big|^{2p}_{L_p(Z,m;\mathcal{L}(H))}ds.
  \end{align*}
  By hypothesis (ii) we have that
  \[
  \big| \mathcal{Q}^h_{1,v(s)} D_1G(u(s),z) -
  D_1^2G(u(s),z)(v(s),\cdot) \big|^{2p}_{\mathcal{L}(H)} \to 0
  \]
  as $h \to 0$ for all $(s,z) \in [0,T]\times Z$, and, by (iii),
  \begin{align*}
    &\E\int_0^T \big| \mathcal{Q}^h_{1,v(s)} D_1G(u(s),\cdot) -
    D_1^2G(u(s),\cdot)(v(s),\cdot)
    \big|^{2p}_{L_p(Z,m;\mathcal{L}(H))} ds\\
    &\qquad \lesssim |h_1|^{2p}_{L_p(Z,m)} \E \int_0^T |v(s)|^{2p}\,ds
    \lesssim |h_1|^{2p}_{L_p(Z,m)} \|v\|_{q'}^p < \infty,
  \end{align*}
  hence, by dominated convergence, $I_{21}(h) \to 0$ as $h \to 0$. In
  a completely similar way one can prove that $I_{22}(h) \to 0$ as $h
  \to 0$ as well. We have thus proved that $I_2(h) \to 0$, hence that
  (\ref{eq:dido}) holds.

  Let us now show that
  \begin{equation}     \label{eq:dudu}
  v \mapsto \partial_{2,v}D_2\F(x,u) \in \L\big(\mathbb{H}_{q'}(T),
  \L(\mathbb{H}_{q'}(T),\mathbb{H}_p(T))\big)
  \end{equation}
  for all $x \in \mathbb{L}_{q'}$ and $u \in \mathbb{H}_{q'}(T)$.  In
  fact, for $w \in \mathbb{H}_{q'}(T)$, we have
  \begin{align*}
    \big\| \partial_{2,v}D_2\F(x,u)w \big\|_p^p &\lesssim
    \E\int_0^T \big| D^2f(u(s))(v(s),w(s))\big|^p\,ds\\\
    &\quad + \E\int_0^T\!\int_Z \big| 
             D_1^2G(u(s),z)(v(s),w(s))\big|^p\,m(dz)\,ds\\
    &\quad + \E\int_0^T\Big(\int_Z \big| 
             D_1^2G(u(s),z)(v(s),w(s))\big|^2\,m(dz)\Big)^{p/2}\,ds\\
    &\lesssim \E\int_0^T |v(s)|^p |w(s)|^p\,ds\\
    &\quad + \big( |h_1|^p_{L_p(Z,m)} + |h_1|^p_{L_2(Z,m)}\big)
       \E\int_0^T |v(s)|^p |w(s)|^p\,ds\\
    &\lesssim \big(1 + |h_1|^p_{L_p(Z,m)} + |h_1|^p_{L_2(Z,m)}\big)
     \|v\|_{q'}^p \|w\|_{q'}^p,
  \end{align*}
  which establishes the continuity of $(v,w)
  \mapsto \partial_{2,v}D_2\F(x,u)w$, and hence ensures that
  (\ref{eq:dudu}) holds true.

  Our next goal is to prove that
  \begin{equation}     \label{eq:fufu}
    u \mapsto \partial_2D_2\F(x,u) \in C\big(\mathbb{H}_{q'}(T),
    \L^{\otimes 2}(\mathbb{H}_{q'}(T),\mathbb{H}_p(T))\big),
  \end{equation}
  for all $x \in \mathbb{L}_{q'}$, which implies the twice continuous
  differentiability of $\F$ by a well-known criterion. Let $u_n \to u$
  in $\mathbb{H}_{q'}(T)$. Then we have
  \begin{align*}
    &\big\| \partial_2 D_2\F(x,u_n)(v,w) - \partial_2 D_2\F(x,u)(v,w)
    \big\|_p^p\\
    &\qquad \lesssim \E\sup_{t\leq T}\Big| \int_0^t e^{(t-s)A}
    [D^2f(u_n(s))(v(s),w(s)) - D^2f(u(s))(v(s),w(s))]\,ds\Big|^p\\
    &\qquad \quad + \E\sup_{t\leq T}\Big| \int_0^t\int_Z e^{(t-s)A}
    [D_1^2G(u_n(s),z)(v(s),w(s)) - D_1^2G(u(s),z)(v(s),w(s))]
    \,\bar{\mu}(ds,dz)\Big|^p\\
    &\qquad =: I_1(n) + I_2(n),
  \end{align*}
  where, using H\"older's inequality with conjugate exponents
  $q'/(2p)$ and $q'/(q'-2p)$,
  \begin{align*}
    I_1(n) &\lesssim \E\int_0^T \big| D^2f(u_n(s))(v(s),w(s))
    - D^2f(u(s))(v(s),w(s))\big|^p\,ds\\
    &\leq \Big(\E\int_0^T \big| D^2f(u_n(s)) - D^2f(u(s))
    \big|_{\L^{\otimes 2}}^{\frac{pq'}{q'-2p}} \Big)^{\frac{q'-2p}{q'}}
    \Big( \E\int_0^T (|v(s)|\,|w(s)|)^{q'/2}\,ds \Big)^{\frac{2p}{q'}}\\
    &\leq \Big(\E\int_0^T \big| D^2f(u_n(s)) - D^2f(u(s))
    \big|_{\L^{\otimes 2}}^{\frac{pq'}{q'-2p}} \Big)^{\frac{q'-2p}{q'}}
    \|v\|_{q'}^p \|w\|_{q'}^p,
  \end{align*}
  which converges to zero as $n \to \infty$ by dominated convergence,
  thanks to the assumption of boundedness of $D^2f$.  Applying again
  inequality (\ref{eq:BDGconv}) yields
  \begin{align*}
  I_2(n) &\lesssim \E\int_0^T\!\int_Z \big| D_1^2G(u_n(s),z)(v(s),w(s))
          - D_1^2G(u(s),z)(v(s),w(s)) \big|^p\,m(dz)\,ds\\
  &\quad + \E\int_0^T\!\Big(\int_Z \big| D_1^2G(u_n(s),z)(v(s),w(s))
          - D_1^2G(u(s),z)(v(s),w(s)) \big|^2\,m(dz)\Big)^{p/2}\,ds\\
  &=: I_{21}(n) + I_{22}(n),
  \end{align*}
  where, again by H\"older's inequality,
  \begin{align*}
  I_{21}(n) &\lesssim \E\int_0^T |v(s)|^p |w(s)|^p \big|
       D_1^2G(u_n(s),\cdot) - D_1^2G(u(s),\cdot)
    \big|^p_{L_p(Z,m;\L^{\otimes 2})}\,ds\\
  &\lesssim \|v\|_{q'}^p \|w\|_{q'}^p \Big( \E\int_0^T\big|
       D_1^2G(u_n(s),\cdot) - D_1^2G(u(s),\cdot)
    \big|^{\frac{pq'}{q'-2p}}_{L_p(Z,m;\L^{\otimes 2})}\,ds \Big)^{\frac{q'-2p}{q'}},
  \end{align*}
  which converges to zero as $n \to \infty$ by continuity of $D_1^2G$
  in its first argument and dominated convergence, thanks to
  hypothesis (iii).  An analogous argument shows that $I_{22}(n) \to
  0$ as $n \to \infty$. We have thus proved (\ref{eq:fufu}). This
  concludes the proof that $\F \in C^2$, hence that the Fr\'echet
  derivative $Du: \mathbb{L}_q \to \L(\mathbb{L}_q,\mathbb{H}_p(T))$
  is G\^ateaux differentiable.

  By Theorem \ref{thm:gattone-impl} we have
  \[
  \partial Du(x)(y_1,y_2) = \big( I - D_2\F(x,u(x)) \big)^{-1}
  D_2^2\F(x,u(x))(Du(x)y_1,Du(x)y_2),
  \]
  hence
  \[
  \partial Du(x)(y_1,y_2) = D_2\F(x,u(x))\partial Du(x)(y_1,y_2) +
  D_2^2\F(x,u(x))(Du(x)y_1,Du(x)y_2),
  \]
  and (\ref{eq:frocione}) now follows substituting in the previous
  identity the expressions for $D_2\F$ and $D_2^2\F$ obtained above
  and in the proof of Theorem \ref{thm:frocetto}.

  The bound for the bilinear form $\partial Du(x)$ can be established
  as an application of Corollary \ref{cor:frocione-impl}. Since
  $Du(x)$ is bounded by Theorem \ref{thm:frocetto}, it is enough to
  show that $D_2^2\F: \mathbb{L}_{q'} \times \mathbb{H}_{q'}(T) \to
  \L^{\otimes 2}(\mathbb{H}_{q'}(T),\mathbb{H}_p(T))$ is bounded. In
  fact, by a computation completely analogous to the above
  ones, we have
  \begin{align*}
  &\big\| D_2^2\F(x,u(x)(v,w) \big\|_p^p \lesssim_{p,T,\eta}
  \int_0^T \big| D^2f(u(s))(v(s),w(s)) \big|^p\,ds\\
  &\qquad + \int_0^T |v(s)|^p |w(s)|^p \Big( 
  \big| D_1^2B(u(s),\cdot) \big|^p_{L_p(Z,m;\L^{\otimes 2})}
  + \big| D_1^2B(u(s),\cdot) \big|^p_{L_2(Z,m;\L^{\otimes 2})}
  \Big)ds\\
  &\quad \lesssim_T C_1^p \|v\|_{q'}^p \|w\|_{q'}^p 
  + \big( |h_1|^p_{L_p(Z,m)} + |h_1|^p_{L_2(Z,m)} \big) \|v\|_{q'}^p \|w\|_{q'}^p.
  \end{align*}

  Let us now assume that $q > 4p \geq 8$, and take $q'\in
  ]2p,q/2[$. We shall deduce the twice Fr\'echet differentiability of
  the solution map applying Theorem \ref{thm:frocione-impl}, with all
  spaces defined as before.  Since $q>2q'$, the first part of the
  proof guarantees that $Du:\mathbb{L}_q \to
  \L(\mathbb{L}_q,\mathbb{H}_{q'}(T))$ is G\^ateaux differentiable
  with derivative $\partial Du: \mathbb{L}_q \to
  \L(\mathbb{L}_q,\L(\mathbb{L}_q,\mathbb{H}_{q'}(T))$.  Therefore, by
  the last statement of Theorem \ref{thm:frocione-impl}, we infer that
  $x \mapsto u(x) \in C^2(\mathbb{L}_q,\mathbb{H}_p(T))$.
\end{proof}


\section{Application: gradient estimates for the resolvent}
\label{sec:ex}
In this section we assume that the coefficients $f$, $B$, and $G$ do
not depend on $t$ and $\omega$. This assumption allows us to
define the semigroup and resolvent associated to the mild solution:
\[
P_t\varphi(x) = \E\varphi(u(t,x)),
\qquad 
R_\alpha\varphi(x) = \int_0^\infty e^{-\alpha t} P_t\varphi(x)\,dt,
\]
where $\varphi\in C_b(H)$ and $\alpha > 0$.

In order to prove gradient estimate for the resolvent $R_\alpha$ we
need the following lemma, which gives an explicit bound on the
Lipschitz constant of the solution map.

\begin{lemma}
Let $A$ be $\eta$-$m$-dissipative, and set
\[
[B]_{1,Q} = \sup_{\substack{x,y \in H\\x\neq y}} \frac{|B(x)-B(y)|_Q}{|x-y|},
\qquad
[G]_{1,m} = \sup_{\substack{x,y \in H\\x\neq y}}
            \frac{|G(x,z)-G(y,z)|_{L_2(Z,m)}}{|x-y|}.
\]
Then we have
\begin{equation}       \label{eq:lippi}
\E|u(t,x) - u(t,y)| \leq e^{\omega_1 t} |x-y|_{\mathbb{L}_2},
\end{equation}
where
\[
\omega_1 := \eta + [f]_1 + \frac12[B]_{1,Q}^2 + \frac12[G]^2_{1,m}
\]
\end{lemma}
\begin{proof}
Let $A_\lambda$ be the Yosida approximation of $A$, $A_\lambda \to A$
as $\lambda \to 0$. Let $u_\lambda$ be the solution of
\[
du_\lambda(t) = (A_\lambda u_\lambda(t) + f(u_\lambda(t)))\,dt
+ B(u_\lambda(t))\,dW(t) + \int_Z G(u_\lambda(t),z)\,\bar{\mu}(dt,dz).
\]
Then It\^o's formula yields
\begin{align*}
|u_\lambda(t)|^2 =& |u_\lambda(0)|^2
+ 2 \int_0^t \ip{A_\lambda u_\lambda(s)}{u_\lambda(s)}\,ds
+ 2 \int_0^t \ip{f(u_\lambda(s))}{u_\lambda(s)}\,ds\\
& + 2 \int_0^t \ip{u_\lambda(s)}{B(u_\lambda(s))dM(s)}
+ \big[B(u_\lambda) \cdot W \big](t)
+ \big[G(u_\lambda) \star \bar{\mu} \big](t)
\end{align*}
hence
\begin{align*}
&\E|u_\lambda(t,x) - u_\lambda(t,y)|^2\\
&\qquad = \E|x-y|^2
+ 2 \E\int_0^t \ip{A_\lambda u_\lambda(s,x)-A_\lambda u_\lambda(s,y)}{%
u_\lambda(s,x) - u_\lambda(s,y)}\,ds\\
&\qquad\quad + 2\E \int_0^t \ip{f(u_\lambda(s,x)) - f(u_\lambda(s,y))}{%
u_\lambda(s,x) - u_\lambda(s,y)}\,ds\\
&\qquad\quad + \E\int_0^t\big|B(u_\lambda(s,x))-B(u_\lambda(s,y))\big|_Q^2\,ds\\
&\qquad\quad + \E\int_0^t \big| G(u_\lambda(s,x),\cdot) - G(u_\lambda(s,y),\cdot)
                    \big|^2_{L_2(Z,m)}ds\\
&\qquad \leq \E|x-y|^2 + \big(2\eta + 2[f]_1 + [B]^2_{1,Q} + [G]^2_{1,m}\big)
\int_0^t \E|u_\lambda(s,x) - u_\lambda(s,y)|^2\,ds
\end{align*}
and Gronwall's inequality implies that
\[
\E|u_\lambda(t,x) - u_\lambda(t,y)|^2 \leq e^{2\omega_1 t} \E|x-y|^2,
\]
therefore, in view of Proposition \ref{prop:yosi}, we can pass to the
limit as $\lambda \to 0$, and applying Cauchy-Schwarz' inequality, we
obtain (\ref{eq:lippi}).
\end{proof}
\begin{rmk}
  Note that (\ref{eq:lippi}) also implies that $|u(x) -
  u(y)|_{\mathcal{H}_2(T)} \leq e^{\omega_1 T} |x-y|_{\mathcal{H}_2}$
  and
  \[
  \E|u(t,x) - u(t,y)| \leq e^{\omega_1 t} |x-y|
  \]
  if $x$, $y\in H$ are nonrandom. Moreover, if the noise is additive,
  i.e. if $B$ is constant, then one can prove, solving differential
  inequalities $\omega$-by-$\omega$, that $|u(t,x) - u(t,y)| \leq
  e^{\omega_1 t} |x-y|$ almost surely.
\end{rmk}
Obtaining an explicit estimate for the Lipschitz constant of the
solution map in the $\mathcal{H}_p(T)$ setting seems considerably more
difficult in the general case of multiplicative noise. It reduces
instead to a simple computation in the case of additive noise, as we
show next.
\begin{lemma}
  Let $A$ be $\eta$-$m$-dissipative and assume that $B$ and $G$ do not
  depend on $x$. If (\ref{eq:musso}) is well-posed in
  $\mathbb{H}_p(T)$, then for any $x$, $y \in \mathcal{L}_p$ we
  have
  \[
  \E|u(t,x) - u(t,y)|^p \leq e^{p\omega_1 t}
  \E|x-y|^p,
  \]
  where $\omega_1=\eta+[f]_1$.
\end{lemma}
\begin{proof}
One has
\[
\frac{d}{dt}\big( u(t,x)-u(t,y) \big) =
Au(t,x)-Au(t,y) + f(u(t,x)) - f(u(t,y))
\]
$\mathbb{P}$-a.s., hence, multiplying both sides by
$|u(t,x)-u(t,y)|^{p-2}(u(t,x)-u(t,y))$, we obtain
\begin{align*}
\frac1p \frac{d}{dt}&\big| u(t,x)-u(t,y) \big|^p\\
&= \big\langle Au(t,x)-Au(t,y),|u(t,x)-u(t,y)|^{p-2}(u(t,x)-u(t,y))
\big\rangle\\
&\quad + |u(t,x)-u(t,y)|^{p-2}
\big\langle f(u(t,x))-f(u(t,y)),(u(t,x)-u(t,y))\big\rangle\\
&\leq \big\langle A|u(t,x)-u(t,y)|^{\frac{p-2}{2}}(u(t,x)-u(t,y)),
|u(t,x)-u(t,y)|^{\frac{p-2}{2}}(u(t,x)-u(t,y)) \big\rangle\\
&\quad + [f]_1 |u(t,x)-u(t,y)|^p\\
&\leq (\eta + [f]_1) |u(t,x)-u(t,y)|^p.
\end{align*}
Writing in integral form and taking expectations, we obtain
\[
\E|u(t,x)-u(t,y)|^p \leq \E|x-y|^p + p(\eta + [f]_1)
\int_0^t \E|u(s,x)-u(s,y)|^p\,ds,
\]
from which the result follows by Gronwall's inequality.
\end{proof}

The following gradient estimate for the resolvent associated to the
mild solution of the stochastic PDE is a consequence of
(\ref{eq:lippi}).
\begin{thm}
  Assume that $\varphi \in C_b(H)$ is G\^ateaux differentiable and
  Lipschitz, and that $f$, $B$, $G$ satisfy the assumptions of Theorem
  \ref{thm:gattino}. Let $\alpha > \omega_1$. Then $x \mapsto
  R_\alpha\varphi(x)$ is G\^ateaux differentiable with
  \begin{equation}
    \label{eq:risolo}
    \partial_y R_\alpha \varphi(x) = \int_0^\infty e^{-\alpha t}
    \E[\partial\varphi(u(t,x))\partial_{2,y}u(t,x)]\,dt,
  \end{equation}
  and it satisfies the estimate
  \[
  \big| \partial R_\alpha\varphi(x) \big|_H \leq
  \frac{1}{\alpha-\omega_1} [\varphi]_1.
  \]
\end{thm}
\begin{proof}
We have, setting $v(t)=\partial_{2,y}u(t,x)$,
\begin{align*}
&\Big| \frac{R_\alpha\varphi(x+hy)-R_\alpha(x)}{h}
- \int_0^\infty e^{-\alpha t} \E[\partial\varphi(u(t,x))\partial_{2,y}u(t,x)]
\,dt \Big|\\
&\qquad \leq \int_0^\infty e^{-\alpha t}
h^{-1}\E[\varphi(u(t,x+hy)) - \varphi(u(t,x))
- h\partial\varphi(u(t,x))v(t)]\,dt\\
&\qquad \leq \int_0^\infty e^{-\alpha t}
\E h^{-1}|\varphi(u(t,x+hy)) - \varphi(u(t,x)+hv(t))|\,dt\\
&\qquad\quad + \int_0^\infty e^{-\alpha t}
\E h^{-1}|\varphi(u(t,x)+hv(t)) - \varphi(u(t,x))
- h\partial\varphi(u(t,x))v(t)|\,dt\\
&\qquad =: I_1(h)+I_2(h),
\end{align*}
and
\[
I_1(h) \leq [\varphi]_1\int_0^\infty e^{-\alpha t}
\E|h^{-1}(u(t,x+hy)-u(t,x))-\partial_{2,y}u(t,x)|\,dt.
\]
Since, by Cauchy-Schwartz' inequality and the differentiability of the
solution map from $H$ to $\mathcal{H}_2(T)$, we have
\begin{align*}
&\E|h^{-1}(u(t,x+hy)-u(t,x))-\partial_{2,y}u(t,x)|\\
&\qquad \leq \big(\E|h^{-1}(u(t,x+hy)-u(t,x))-\partial_{2,y}u(t,x)|^2\big)^{1/2}
\to 0,
\end{align*}
we conclude that $I_1(h) \to 0$ as $h \to 0$ by dominated convergence.
On the other hand, $I_2(h)$ converges to zero as $h \to 0$ by
definition of directional derivative and dominated convergence. This
establishes (\ref{eq:risolo}).

Note that $x \mapsto P_t\varphi(x)$ is Lipschitz: in fact, the
previous lemma yields
\begin{align*}
\big| P_t\varphi(x) - P_t\varphi(y) \big| &\leq
\E\big| \varphi(u(t,x)) - \varphi(u(t,y)) \big|\\
&\leq [\varphi]_1 \E|u(t,x)-u(t,y)| \leq [\varphi]_1 e^{\omega_1t} |x-y|.
\end{align*}
Therefore
\begin{align*}
\frac{R_\alpha\varphi(x+hy)-R_\alpha\varphi(x)}{h} &=
\int_0^\infty e^{-\alpha t} h^{-1}
\big(P_t\varphi(x+hy)-P_t\varphi(x)\big)\,dt\\
&\leq [\varphi]_1 |y| \int_0^\infty  e^{\omega_1t} e^{-\alpha t} \,dt < \infty,
\end{align*}
hence, in view of (\ref{eq:risolo}),
\[
\big| \partial_y R_\alpha\varphi(x) \big| \leq 
\frac{1}{\alpha-\omega_1} [\varphi]_1 |y|,
\]
thus also $|\partial R_\alpha\varphi(x)| \leq
(\alpha-\omega_1)^{-1}[\varphi]_1$.
\end{proof}


\section{Strong Feller property}     \label{sec:Feller}
The purpose of this section is to establish a Bismut-Elworthy formula
for the semigroup associated to (\ref{eq:musso}), and to deduce from
it the strong Feller property, adapting an argument of \cite{PriZa} to
the infinite dimensional case. We would like to emphasize that the
proof depends essentially on the second order differentiability of the
solution with respect to the initial datum established in Theorem
\ref{thm:frocione} above. In the following we shall denote the set of
bounded Borel functions from $H$ to $\erre$ by $B_b(H)$.
\begin{thm}
  Assume that $Q \in \mathcal{L}_2(K)$, $B(x)$ is invertible with
  $|B(x)^{-1}| \leq C$ for all $x\in H$, for some $C>0$, and the
  hypotheses of Theorem \ref{thm:gattino} are satisfied. Then the
  semigroup $P_t$ is strong Feller, i.e. $\varphi \in
  B_b(H)$ implies $P_t\varphi \in C_b(H)$.
\end{thm}
\begin{proof}
  We first assume that the coefficients $f$, $B$ and $G$ satisfy the
  hypotheses of Theorem \ref{thm:frocione}, so that $x \mapsto u(x)
  \in C^2(H,\mathbb{H}_2(T))$. This assumption will be removed in the
  last part of the proof.

  A formal application of It\^o's formula shows that the generator $L$
  of the semigroup $P_t$ associated to the mild solution of
  (\ref{eq:musso}) takes the form, for $\varphi\in C^2_b(H)$,
  \begin{align*}
    L\varphi(x) =& \ip{Ax+f(x)}{D\varphi(x)}
    + \frac12 \tr(QB(x)B^*(x)D^2\varphi(x))\\
    & + \int_Z \big[ \varphi(x+G(x,z)) - \varphi(x) -
    \ip{D\varphi(x)}{G(x,z)} \big]\,m(dz).
  \end{align*}
  Let $u_\lambda$ be the solution of (\ref{eq:musso}) with $A$
  replaced by its Yosida approximation $A_\lambda$, $P_t^\lambda$ the
  associated semigroup, and $L_\lambda$ the generator of
  $P_t^\lambda$. Then the action of $L_\lambda$ on $\varphi\in
  C^2_b(H)$ is exactly as for $L$, with $A$ replaced by $A_\lambda$.
  Let $\varphi \in C^2_b(H)$, $s \in [0,t]$, and set
  $v(s,x)=P^\lambda_{t-s}\varphi(x)\equiv
  \E\varphi(u_\lambda(t-s,x))$. Then $v \in C^{1,2}([0,T]\times H)$
  and It\^o's formula implies
  \begin{align*}
    v(s,u_\lambda(s)) =\;& v(0,x) 
    + \int_0^s (\partial_r + L_\lambda)v(r,u_\lambda(r))\,dr
    + \int_0^s \ip{Dv(r,u_\lambda(r))}{B(u_\lambda(r))\,dW(r)}\\
    & + \int_0^s\!\int_Z \big[v(r-,u_\lambda(r-)+G(u_\lambda(r-),z)) -
    v(r-,u_\lambda(r-))\big]\, \bar{\mu}(dr,dz)
  \end{align*}
  Since $(\partial_t + L_\lambda)v=0$, the previous identity evaluated at
  $s=t$ implies
  \[
    \varphi(u_\lambda(t)) = P^\lambda_t\varphi(x) + M^\lambda_1(t)
    + M^\lambda_2(t),
  \]
  where
  \begin{eqnarray*}
    M^\lambda_1(t) &=& \int_0^t \ip{DP^\lambda_{t-r}\varphi(u_\lambda(r))}{%
    B(u_\lambda(r))dW(r)},\\
    M^\lambda_2(t) &=& \int_0^t\!\int_Z 
    \big[ P^\lambda_{t-s}\varphi(u_\lambda(s-)+G(u_\lambda(s-))z)
    - P^\lambda_{t-s}\varphi(u_\lambda(s-)) \big]
    \,\bar\mu(ds,dz).
  \end{eqnarray*}
  Letting $\lambda \to 0$ and recalling Proposition \ref{prop:yosi} we
  obtain
  \begin{equation}     \label{eq:deco}
  \varphi(u(t)) = P_t\varphi(x) + M_1(t) + M_2(t),
  \end{equation}
  with $M_1$ and $M_2$ defined in the obvious way.
  Moreover, setting $w(t)=\partial_{2,y}u(t,x)$ and
  \[
  M_3(t) = \int_0^t \ip{B^{-1}(u(s))w(s)}{dW(s)},
  \]
  multiplying both sides of (\ref{eq:deco}) by $M_3(t)$ and taking
  expectations yields
  \begin{eqnarray*}
    \E\varphi(u(t))M_3(t) &=& \E M_1(t)M_3(t) =
    \E\int_0^t \ip{DP_{t-s}\varphi(u(s))}{w(s)}\,ds\\
    &=& \E\int_0^t D[P_{t-s}\varphi(u(s))]y\,ds =
    \int_0^t DP_t\varphi(x)y\,ds = t DP_t\varphi(x)y.
  \end{eqnarray*}
  Here $\E M_2(t)M_3(t)=0$ because $W$ and $\bar\mu$ are independent,
  and we have used the Markov property of solutions in the second to last
  step.  In particular, we have proved the Bismut-Elworthy-type
  formula
  \[
  DP_t\varphi(x)y = \frac1t \E\Big[ \varphi(u(t,x)) \int_0^t
  \ip{B^{-1}(u(s,x))\partial_{2,y}u(s,x)}{dW(s)} \Big].
  \]
  We shall now remove the assumptions on $f$, $B$ and $G$. Let us
  assume for a moment that we can find sequences $f_\varepsilon$,
  $B_\varepsilon$, $G_\varepsilon$ satisfying the hypotheses of
  Theorem \ref{thm:frocione} and Proposition \ref{prop:vava}, and denote
  the mild solution of
  \[
  du(t) = [Au(t) + f_\varepsilon(u(t))]\,dt +
  B_\varepsilon(u(t))\,dW(t) + \int_Z
  G_\varepsilon(u(t),z)\,\bar{\mu}(dt,dz), \qquad u(0)=x,
  \]
  by $u_\varepsilon$, so that $x \mapsto u_\varepsilon(x) \in
  C^2(H,\mathcal{H}_2(T))$ and $u_\varepsilon \to u$ in
  $\mathcal{H}_2(T)$. In particular we also have
  $P^\varepsilon_t\varphi(x) \to P_t\varphi(x)$ for all $x\in H$ and
  $t \leq T$, where
  $P^\varepsilon_t\varphi(x):=\E\varphi(u_\varepsilon(t,x))$, $\varphi
  \in C_b(H)$. Then Cauchy-Schwartz' inequality yields
  \[
  |DP^\varepsilon_t\varphi(x)y|^2 \leq \frac{1}{t^2}
  |\varphi|^2_\infty C^2 \E \int_0^t |\partial_{2,y}u_\varepsilon(s,x)|^2\,ds,
  \]
  where $|\varphi|_\infty:=\sup_{x\in H} |\varphi(x)|$.  In view of
  Remark \ref{rmk:lipco} it is not difficult to see that there exists
  a constant $N$, which does not depend on $x$, $y$, and
  $\varepsilon$, such that $|[u_\varepsilon(x_1)-u_\varepsilon(x_2]|_2
  \leq N|x_1-x_2|_H$, hence, by Theorem \ref{thm:gattino},
  $|\partial_{2,y}u_\varepsilon(s,x)| \leq N|y|$. We obtain
  \[
  |DP^\varepsilon_t\varphi(x)| \leq \frac{NC}{t^{1/2}}
  |\varphi|_{\infty},
  \]
  thus also
  $|P^\varepsilon_t\varphi(x_1)-P^\varepsilon_t\varphi(x_2)| \leq
  t^{-1/2}NC|\varphi|_{\infty}|x_1-x_2|$, and letting $\varepsilon \to
  0$,
  \[
  |P_t\varphi(x_1)-P_t\varphi(x_2)| \leq
  t^{-1/2}NC|\varphi|_{\infty}|x_1-x_2|.
  \]
  The same Lipschitz property continues to hold also for $\varphi \in
  B_b(H)$ by a simple regularization argument
  (see e.g. \cite[Lemma 2.2]{PZ95}).

  In order to complete the proof, we have to show that we can find
  sequences $f_\varepsilon$, $B_\varepsilon$, $G_\varepsilon$
  satisfying the hypotheses of Theorem \ref{thm:frocione} and
  Proposition \ref{prop:vava}. The existence of such $f_\varepsilon$
  and $B_\varepsilon$ is well-known (see
  e.g. \cite[Sect.~3.3.1]{DP-K}, \cite{PZ95}), and the construction of
  $G_\varepsilon$ can be carried out in a completely similar way, hence
  we omit it.
\end{proof}

\subsection*{Acknowledgments}
The first author was partially supported by the SFB 611, Bonn, by the
ESF through grant AMaMeF 969, by the Centre de Recerca Matem\`atica,
Barcelona, through an EPDI fellowship, and by the EU through grant
MOIF-CT-2006-040743. The two last named authors were supported by SFB
701, NSF Grant 0606615 and the BiBoS Research Centre in Bielefeld.

\let\oldbibliography\thebibliography
\renewcommand{\thebibliography}[1]{%
  \oldbibliography{#1}%
  \setlength{\itemsep}{-1pt}%
}

\bibliographystyle{amsplain}
\bibliography{ref}

\parindent=0pt

Carlo Marinelli, Institut f\"ur Angewandte Mathematik, Universit\"at Bonn,
Wegelerstr. 6, D-53115 Bonn, Germany, and Facolt\`a di Economia,
Libera Universit\`a di Bolzano, Via Sernesi 1, I-39100 Bolzano.\\
e-mail \verb+cm788@uni-bonn.de+

\medskip

Claudia Pr\'ev\^ot, Fakult\"at f\"ur Mathematik, Universit\"at
Bielefeld, Postfach 100 131, D-33501 Bielefeld, Germany.

\medskip

Michael R\"ockner, Fakult\"at f\"ur Mathematik, Universit\"at
Bielefeld, Postfach 100 131, D-33501 Bielefeld, Germany, and
Departments of Mathematics and Statistics, Purdue University,
150 N. University St., West Lafayette, IN 47907-2067, USA.\\
e-mail \verb+roeckner@math.uni-bielefeld.de+

\end{document}